\documentclass[11pt]{article}

\usepackage[T1]{fontenc}
\usepackage[utf8]{inputenc}
\usepackage{lmodern}
\usepackage{microtype}
\usepackage{indentfirst}
\usepackage[a4paper,margin=1in]{geometry}
\usepackage{amsmath,amssymb,amsthm,mathtools}
\usepackage{booktabs}
\usepackage{xcolor}
\usepackage{hyperref}
\usepackage[nameinlink,capitalise,noabbrev]{cleveref}

\hypersetup{colorlinks=true,linkcolor=blue,citecolor=blue,urlcolor=blue,
  pdfauthor={Zhekai Pang}}
\makeatletter
\renewcommand{\tagform@}[1]{%
  \maketag@@@{\color{blue}(\ignorespaces#1\unskip\@@italiccorr)}%
}
\makeatother

\newtheorem{theorem}{Theorem}[section]
\newtheorem{lemma}[theorem]{Lemma}
\newtheorem{proposition}[theorem]{Proposition}
\newtheorem{corollary}[theorem]{Corollary}
\theoremstyle{definition}

\theoremstyle{remark}

\newcommand{\dd}{\,\mathrm{d}}

\title{Entropy maximization for cube projections}
\hypersetup{pdftitle={Entropy maximization for cube projections}}
\author{Zhekai Pang\\
\small Universitat Pompeu Fabra, Barcelona, Spain\\
\small\texttt{zhekai.pang@upf.edu}}
\date{}
\begin{document}
\maketitle
\begin{abstract}
We prove that, in every dimension $n\ge150$, the one-dimensional
projections of the uniform distribution on a cube have the largest
differential entropy in directions with equal coordinate magnitudes.
This establishes the conjecture of Eskenazis, Nayar, and Tkocz in
these dimensions.
\end{abstract}



\section{Introduction}

Let $U_1,\ldots,U_n$ be independent and uniform on $[-1,1]$.
The one-dimensional projections of the uniform distribution
on the cube $[-1,1]^n$ are the sums $\sum\limits_{i=1}^n a_iU_i$.
Under the normalization $\sum\limits_{i=1}^n a_i^2=1$, these sums all have
variance $1/3$. Their entropies depend on the coefficients,
and we ask which choices maximize the entropy.

Eskenazis, Nayar, and Tkocz \cite[Question~9]{ENT2018}
conjectured that the maximum is attained when all coefficient
magnitudes are equal. The same question appears in Bartczak,
Nayar, and Zwara \cite[Question~3]{BNZ2021} and Chasapis,
Gurushankar, and Tkocz \cite[Section~6]{CGT2023}.
We prove this conjecture for $n\ge150$.

\begin{theorem}\label{thm:q3-main}
Let $n\ge150$ be an integer, and let $U_1,\ldots,U_n$ be independent
random variables, each uniform on $[-1,1]$. For a random variable $X$
with density $f$, let
\[
 h(X)=-\int_{-\infty}^{\infty}f(x)\log f(x)\dd x.
\]
For every $a\in\mathbb R^n$ with $\sum\limits_{i=1}^n a_i^2=1$,
\[
 h\left(\sum\limits_{i=1}^n a_iU_i\right)
 \le h\left(\frac1{\sqrt n}\sum\limits_{i=1}^n U_i\right).
\]
Equality holds if and only if $|a_i|=1/\sqrt n$ for every $i$.
\end{theorem}

The same maximization problem has been solved for several
other families of summands. Eskenazis, Nayar, and Tkocz
\cite[Theorem~8]{ENT2018} proved that equal coefficient magnitudes
maximize the entropy of sums of independent copies of a Gaussian
scale mixture, provided the entropies are finite. Bartczak, Nayar,
and Zwara \cite{BNZ2021} proved the corresponding result for
nonnegative coefficients and independent $\operatorname{Gamma}(\alpha)$
summands when $n\alpha\ge1$.

The projection density is even by symmetry and log-concave by
Pr\'ekopa's marginal theorem \cite{Prekopa1973}, so it attains its
maximum at zero. For a unit projection direction, the coarea
formula identifies this value with $2^{-n}$ times the
$(n-1)$-dimensional volume of the central section perpendicular
to that direction. Ball's cube slicing theorem \cite{Ball1986}
maximizes this value over unit directions, with a maximizer having
exactly two nonzero coordinates of equal magnitude.
Melbourne and Roberto \cite{MR2022} proved a quantitative version
and characterized equality in the associated min-entropy power
inequality. For Shannon entropy, Ball, Nayar, and Tkocz
\cite{BNT2016} proved a reverse entropy power inequality for
one-dimensional projections of general log-concave random vectors.

Theorem~\ref{thm:q3-main} compares the Shannon entropies of the full
projection densities as the weights vary in a fixed dimension.
A separate question is how the entropy of the equally weighted
sum changes with $n$. The result of Artstein, Ball, Barthe, and Naor
\cite{ABBN2004}, applied to uniform summands, shows that
$h(n^{-1/2}\sum\limits_{i=1}^n U_i)$ is nondecreasing in $n$.
Barron \cite{Barron1986} proved convergence to Gaussian entropy,
and Bobkov, Chistyakov, and G\"otze \cite{BCG2013} obtained
Edgeworth expansions for the relative entropy in this limit.
Brazitikos, Rapaport, and Tkocz \cite{BRT2026} gave cumulant
bounds for relative Fisher information and Gaussian entropy
deficits, with applications to weighted uniform sums.

One approach is to average pairs of squared coefficients.
Lata\l a and Oleszkiewicz \cite{LO1995} used this method
to obtain sharp moment inequalities for uniform sums.
Chasapis, Gurushankar, and Tkocz \cite{CGT2023}
proved sharp lower bounds on the $p$-norm for $0<p<1$.
Malliaris \cite{Malliaris2025} gave Fourier criteria for
monotonicity of expectations under coefficient balancing.
The principal criterion in that work assumes that the Fourier
transform of the summand density is positive and integrable.
The uniform transform $\sin t/t$ changes sign, so this criterion
does not directly give the entropy comparison considered here.

Let $q_n$ be the density of $n^{-1/2}\sum\limits_{i=1}^n U_i$.
Bartczak, Nayar, and Zwara \cite[Section~8.3.1]{BNZ2021}
observed that the entropy comparison follows if
$(\log q_n)'''(x)\le0$ for every $x\in(0,\sqrt n)$.
For $n\ge150$ and every $x\in(0,\sqrt n)$, we prove the strict inequality
\[
 (\log q_n)'''(x)<0.
\]

Exponential tilting associates each point of $(0,\sqrt n)$ with
a positive tilting parameter and reduces the sign of
$(\log q_n)'''$ there to the sign of the third logarithmic derivative
at the origin for a centered and standardized sum. We use real Fourier inversion to estimate
the density derivatives at the origin for every positive tilting
parameter. As the tilt tends to zero, the third moment of the
centered summand also tends to zero. The centered tilted summand
density minus its reflection changes sign once on the positive
half-line. We use this sign change to obtain a comparison that
retains the magnitude of the third moment as a factor in the error
bounds for the odd derivatives.

Section~\ref{sec:cube-comparison} gives the entropy reduction.
Sections~\ref{sec:q3-tilt-reduction}--\ref{sec:q3-global-odd-tail}
prove the derivative inequality.
Section~\ref{sec:information-theory} provides some applications to information theory.

\section{Entropy comparison}
\label{sec:cube-comparison}

We first prove a comparison for two coefficients, using the
averaging argument of Lata\l a and Oleszkiewicz \cite{LO1995}.

\begin{lemma}\label{lem:cube-balancing}
Let $U_1,U_2$ be independent random variables, each uniform on
$[-1,1]$. Let $V$ be even and $C^3$ on $(-R,R)$, with
$V'''(x)\ge0$ for $0<x<R$. Let $Y$ be independent of $U_1,U_2$,
and suppose that either $Y=0$ almost surely or $Y$ has an even
unimodal density.
If $a\ge b\ge0$, $c=\sqrt{(a^2+b^2)/2}$, and $|Y|\le L$ almost
surely with $2c+L<R$, then
\[
 \mathbb E[V(aU_1+bU_2+Y)]
 \le\mathbb E[V(cU_1+cU_2+Y)].
\]
\end{lemma}

\begin{proof}
The assertion is an identity when $a=b$, so suppose that $a>b$.
Let $g(x)=\mathbb E[V(x+Y)]$ for $|x|\le2c$.
Because $2c+L<R$, the arguments $x+Y$ stay in
$[-2c-L,2c+L]\subseteq(-R,R)$. The first three derivatives of $V$
are bounded on a slightly larger interval, so differentiation
under the expectation is valid through order three.
If $Y=0$ almost surely, then $g=V$ and $g'''\ge0$ on $[0,2c]$.

Otherwise, let $q$ be the even unimodal density of $Y$.
For $x,u\ge0$, the inequality $|u-x|\le u+x$ gives
$q(u-x)\ge q(u+x)$. To compute $g'''$, we first change variables
by $u=x+y$. Since $[x-L,x+L]\subseteq(-R,R)$ for $x\in[0,2c]$,
the support of $q$ permits integration over $(-R,R)$.
Splitting this integral at zero and using the oddness of
$V'''$ and the evenness of $q$, for every $x\in[0,2c]$ we obtain
{
\medmuskip=1mu\relax
\thickmuskip=2.75mu\relax
\begin{align*}
 g'''(x)
 &=\int_{-L}^{L}\!V'''(x+y)q(y)\mkern1mu\mathrm{d}y
 =\int_{-R}^{R}\!V'''(u)q(u-x)\mkern1mu\mathrm{d}u
 =\int_0^R\!V'''(u)[q(u-x)-q(u+x)]\mkern1mu\mathrm{d}u\ge0.
\end{align*}
}
In either case, $g''$ is nondecreasing on $[0,2c]$.
Since $V$ is even and $Y$ is symmetric, $g$ is even and
$g'(0)=0$. Writing
$g'(x)=\int_0^x g''(v)\,\mathrm{d}v$, for every $x\in[0,2c]$
we obtain
\begin{equation}\label{eq:cube-square-convexity}
 xg''(x)-g'(x)=\int_0^x[g''(x)-g''(v)]\,\mathrm{d}v\ge0.
\end{equation}
Let $\psi(t)=g(\sqrt t)$ for $t\in[0,4c^2]$. For every
$t\in(0,4c^2)$, \eqref{eq:cube-square-convexity} gives
\[
 \psi''(t)=\frac{\sqrt t\,g''(\sqrt t)-g'(\sqrt t)}{4t^{\frac32}}\ge0.
\]
Continuity at $0$ and $4c^2$ therefore makes $\psi$ convex on
$[0,4c^2]$.

Suppose first that $a>b>0$. Let $f_1$ and $f_2$ be the densities of
$aU_1+bU_2$ and $c(U_1+U_2)$, respectively. Both densities are
even. Since $a+b<2c$, they vanish outside $[-2c,2c]$.
By symmetry, $\int_0^{2c}f_j(x)\,\mathrm{d}x=1/2$ for $j=1,2$.
Since $\mathbb E[U_i]=0$ and $\mathbb E[U_i^2]=1/3$, independence gives
\[
 \int_0^{2c}x^2f_1(x)\,\mathrm{d}x
 =\frac{a^2+b^2}{6}=\frac{c^2}{3}
 =\int_0^{2c}x^2f_2(x)\,\mathrm{d}x.
\]
Subtracting the two masses and the two second moments yields
\begin{equation}\label{eq:cube-pair-moments}
 \int_0^{2c}(f_2-f_1)(x)\,\mathrm{d}x
 =\int_0^{2c}x^2(f_2-f_1)(x)\,\mathrm{d}x=0.
\end{equation}
On $[0,2c]$, their density functions are
\begin{equation}\label{eq:cube-pair-densities}
 \begin{aligned}
 f_1(x)&=
 \begin{cases}
 (2a)^{-1},&0\le x\le a-b,\\
 \frac{a+b-x}{4ab},&a-b<x<a+b,\\
 0,&a+b\le x\le2c,
 \end{cases}\\
 f_2(x)&=\frac{2c-x}{4c^2}.
 \end{aligned}
\end{equation}
Since $c<a$ and $c^2>ab$, \eqref{eq:cube-pair-densities} shows that
$f_2-f_1$ starts positive, strictly decreases on $[0,a-b]$, and
strictly increases on $[a-b,a+b]$. On $[a+b,2c]$ it is nonnegative.
The zero integral in \eqref{eq:cube-pair-moments} forces $f_2-f_1$
to be negative somewhere. It therefore changes sign exactly twice
in $(0,2c)$, with the sign pattern $+,-,+$.
Let $r_1<r_2$ be these two sign-change points.

Let $\ell$ be the affine function interpolating $\psi$ at
$r_1^2,r_2^2$. Convexity gives $\psi(x^2)\le\ell(x^2)$ for
$x\in[r_1,r_2]$ and $\psi(x^2)\ge\ell(x^2)$ for
$x\in[0,r_1]\cup[r_2,2c]$.
Thus the product of $\psi(x^2)-\ell(x^2)$ and $(f_2-f_1)(x)$
is nonnegative.
Because $\ell$ is affine, $\ell(x^2)$ is a constant plus a
multiple of $x^2$. The two equalities in
\eqref{eq:cube-pair-moments} therefore give
\[
 \int_0^{2c}\ell(x^2)(f_2-f_1)(x)\,\mathrm{d}x=0.
\]
By the definition of $\psi$, we have $g(x)=\psi(x^2)$ for
$x\in[0,2c]$. Substituting this identity and using
\eqref{eq:cube-pair-moments}, we obtain
\begin{equation}
 \int_0^{2c}g(x)(f_2-f_1)(x)\,\mathrm{d}x
 =\int_0^{2c}\psi(x^2)(f_2-f_1)(x)\,\mathrm{d}x
 =\int_0^{2c}[\psi(x^2)-\ell(x^2)](f_2-f_1)(x)\,\mathrm{d}x\ge0.
 \label{eq:cube-pair-comparison}
\end{equation}
Since $Y$ is independent of $U_1,U_2$, the definition of $g$ gives
\[
 \mathbb E[V(aU_1+bU_2+Y)]=\mathbb E[g(aU_1+bU_2)],\qquad
 \mathbb E[V(cU_1+cU_2+Y)]=\mathbb E[g(cU_1+cU_2)].
\]
The functions $g,f_1,f_2$ are even, and both densities vanish
outside $[-2c,2c]$. Hence \eqref{eq:cube-pair-comparison} yields
\[
 \mathbb E[g(cU_1+cU_2)]-\mathbb E[g(aU_1+bU_2)]
 =\int_{-2c}^{2c}g(x)(f_2-f_1)(x)\,\mathrm{d}x
 =2\int_0^{2c}g(x)(f_2-f_1)(x)\,\mathrm{d}x\ge0.
\]
Combining these identities proves the claimed inequality when $a>b>0$.
When $b=0<a$, apply this case with coefficients $a$ and
$\varepsilon>0$, whose balanced coefficient is
$\sqrt{(a^2+\varepsilon^2)/2}$. For sufficiently small $\varepsilon$,
the strict support bound still holds. All arguments of $V$
then stay in a fixed compact subinterval of $(-R,R)$.
Continuity and the dominated convergence theorem allow
$\varepsilon\downarrow0$ in both expectations.
\end{proof}

We now use Lemma~\ref{lem:cube-balancing} to prove the entropy
criterion proposed by Bartczak, Nayar, and Zwara
\cite[Section~8.3.1]{BNZ2021}.

\begin{proposition}\label{prop:cube-criterion}
Let $n\ge150$, and let $q_n$ be the density of
$n^{-1/2}\sum\limits_{i=1}^n U_i$. Suppose that, for every
$x\in(0,\sqrt n)$,
\begin{equation}\label{eq:cube-third}
 (\log q_n)'''(x)\le0.
\end{equation}
Then, for every $a\in\mathbb R^n$ with $\sum\limits_{i=1}^n a_i^2=1$,
\[
 h\left(\sum\limits_{i=1}^n a_iU_i\right)
 \le h\left(\frac1{\sqrt n}\sum\limits_{i=1}^n U_i\right),
\]
with equality if and only if $|a_i|=1/\sqrt n$ for every $i$.
\end{proposition}

\begin{proof}
Since each $U_i$ is symmetric about zero, changing the signs of
the coefficients does not change the distribution of the sum.
We may therefore assume that $a_i\ge0$ for every $i$.
Let $V=-\log q_n$.
We compare expectations of $V$ by averaging coefficients,
then use nonnegativity of relative entropy to obtain the
entropy bound.

\medskip
\noindent\emph{Averaging the coefficients.}
Convolution of the uniform densities shows that $q_n$ is even,
positive on $(-\sqrt n,\sqrt n)$, and has continuous derivatives
through order $n-2$. Since $n\ge150$, the function $V$ is even
and $C^3$ on this interval. The hypothesis
\eqref{eq:cube-third} gives $V'''(x)\ge0$ for $0<x<\sqrt n$.
The function $V$ tends to $+\infty$ at both endpoints, so we
cannot apply Lemma~\ref{lem:cube-balancing} directly when the
support reaches $\pm\sqrt n$. To obtain the strict support
bound required by the lemma, fix $0<\rho<1$ and consider
$x\mapsto V(\rho x)$. This function is even and $C^3$ on
$(-\sqrt n/\rho,\sqrt n/\rho)$, has nonnegative third derivative
on $(0,\sqrt n/\rho)$, and is bounded on $[-\sqrt n,\sqrt n]$.

For this fixed $\rho$, we first show that replacing any two coefficients by their
common root mean square does not decrease
\[
 \mathbb E\left[V\left(\rho\sum\limits_{i=1}^n a_iU_i\right)\right].
\]
Choose any two coefficients and keep all others fixed.
Since the $U_i$ are independent and identically distributed,
we may relabel the coefficients and the corresponding variables
so that the chosen pair is $a_1\ge a_2$.
We replace $a_1,a_2$ by $c,c$, where $c=\sqrt{(a_1^2+a_2^2)/2}$.
Thus the two squared coefficients are replaced by their
arithmetic mean. To apply Lemma~\ref{lem:cube-balancing} to
this replacement, let
\[
 Y=\sum\limits_{i=3}^n a_iU_i.
\]

$Y$ is independent of $U_1,U_2$. If $a_i=0$ for every $3\le i\le n$, then $Y=0$
almost surely. Otherwise, we obtain its density by adding the
nonzero summands one at a time. The first such summand is
uniform on an interval centered at zero, so its density is
even and nonincreasing on $[0,\infty)$.

Suppose that a partial sum has a density $q$ with these
properties, and let $b>0$ be the coefficient of the next
nonzero summand. Independence gives the density of the new
partial sum as
\[
 x\longmapsto\frac1{2b}\int_{x-b}^{x+b}q(u)\,\mathrm{d}u.
\]
This density is even because $q$ is even. For $x>0$, the
inequality $|x-b|\le x+b$ and the monotonicity of $q$ give
$q(x-b)=q(|x-b|)\ge q(x+b)$. Hence, for almost every $x>0$,
\[
 \frac{\mathrm{d}}{\mathrm{d}x}
 \left[\frac1{2b}\int_{x-b}^{x+b}q(u)\,\mathrm{d}u\right]
 =\frac{q(x+b)-q(x-b)}{2b}\le0.
\]
The new density is absolutely continuous, so it is also
nonincreasing on $[0,\infty)$. Repeating this argument shows
that $Y$ has an even unimodal density.
The same convolution argument shows that $q_n$ is even and
nonincreasing on $[0,\infty)$.

Since $|U_i|\le1$ and $a_i\ge0$, we have
$|Y|\le\sum\limits_{i=3}^n a_i$ almost surely.
The definition of $c$ gives
$2c^2+\sum\limits_{i=3}^n a_i^2=\sum\limits_{i=1}^n a_i^2=1$.
The Cauchy--Schwarz inequality applied to the new coefficient
vector therefore yields
\[
 2c+\sum\limits_{i=3}^n a_i\le\sqrt n<\frac{\sqrt n}{\rho}.
\]
All assumptions of Lemma~\ref{lem:cube-balancing} now hold for
$x\mapsto V(\rho x)$ and $Y$.
The lemma gives
\begin{equation}\label{eq:cube-single-average}
 \mathbb E[V(\rho(a_1U_1+a_2U_2+Y))]
 \le\mathbb E[V(\rho(cU_1+cU_2+Y))].
\end{equation}

Starting with $a_i^{(0)}=a_i$, obtain each new coefficient
vector by replacing a largest and a smallest coefficient
by their common root mean square.
All coefficients remain nonnegative, and
$\sum\limits_{i=1}^n(a_i^{(k)})^2=1$ at every step.
For every integer $k\ge0$, repeated application of
\eqref{eq:cube-single-average} gives
\[
 \mathbb E\left[V\left(\rho\sum\limits_{i=1}^n a_iU_i\right)\right]
 \le\mathbb E\left[V\left(\rho\sum\limits_{i=1}^n a_i^{(k)}U_i\right)\right].
\]

We show that the squared coefficients converge to their
common mean $1/n$ by computing how much their sum of squared
deviations decreases at each step.
Fix an integer $k\ge0$, and let
\[
 x=\max\limits_{1\le i\le n}(a_i^{(k)})^2,
 \qquad y=\min_{1\le i\le n}(a_i^{(k)})^2.
\]
At the next step, these two squared coefficients are both
replaced by $(x+y)/2$. All other squared coefficients are
unchanged, so their contributions cancel when we subtract
the two sums, giving
{
\medmuskip=2.5mu\relax
\thickmuskip=3.5mu\relax
\begin{align*}
 \sum\limits_{i=1}^n\left[(a_i^{(k)})^2-\frac1n\right]^2
 -\sum\limits_{i=1}^n\left[(a_i^{(k+1)})^2-\frac1n\right]^2
 &=\left(x-\frac1n\right)^2
          +\left(y-\frac1n\right)^2
          -2\left(\frac{x+y}{2}-\frac1n\right)^2
 =\frac12(x-y)^2.
\end{align*}
}

The mean $1/n$ and all the squared coefficients lie in
$[y,x]$. Thus, for every $1\le i\le n$,
\[
 \left|(a_i^{(k)})^2-\frac1n\right|\le x-y.
\]
Squaring and summing gives
\[
 \sum\limits_{i=1}^n\left[(a_i^{(k)})^2-\frac1n\right]^2
 \le n(x-y)^2.
\]
Hence
\[
 \sum\limits_{i=1}^n\left[(a_i^{(k+1)})^2-\frac1n\right]^2
 =\sum\limits_{i=1}^n\left[(a_i^{(k)})^2-\frac1n\right]^2
   -\frac12(x-y)^2
 \le\left(1-\frac1{2n}\right)
   \sum\limits_{i=1}^n\left[(a_i^{(k)})^2-\frac1n\right]^2.
\]
Iterating this inequality, as $k\to\infty$ we obtain
\[
 0\le\sum\limits_{i=1}^n\left[(a_i^{(k)})^2-\frac1n\right]^2
 \le\left(1-\frac1{2n}\right)^k
     \sum\limits_{i=1}^n\left(a_i^2-\frac1n\right)^2
 \longrightarrow0.
\]
Since the coefficients are nonnegative, this implies
$a_i^{(k)}\to1/\sqrt n$ for every $1\le i\le n$.
Using the same $U_i$ at every step, we obtain
\[
 \sum\limits_{i=1}^n a_i^{(k)}U_i
 \longrightarrow\frac1{\sqrt n}\sum\limits_{i=1}^n U_i
 \quad\text{almost surely}.
\]
Since $\sum\limits_{i=1}^n(a_i^{(k)})^2=1$ and $|U_i|\le1$ almost surely,
for every integer $k\ge0$ the Cauchy--Schwarz inequality gives
\[
 \left|\sum\limits_{i=1}^n a_i^{(k)}U_i\right|
 \le\sqrt{\sum\limits_{i=1}^n(a_i^{(k)})^2}
     \sqrt{\sum\limits_{i=1}^n U_i^2}
 \le\sqrt n.
\]
The function $x\mapsto V(\rho x)$ is bounded and
continuous on $[-\sqrt n,\sqrt n]$, so the dominated convergence
theorem allows us to pass to the limit after iterating
\eqref{eq:cube-single-average}, giving
\begin{equation}\label{eq:cube-contracted-comparison}
 \mathbb E\left[V\left(\rho\sum\limits_{i=1}^n a_iU_i\right)\right]
 \le\mathbb E\left[V\left(\frac{\rho}{\sqrt n}\sum\limits_{i=1}^n U_i\right)\right].
\end{equation}

\medskip
\noindent\emph{Letting $\rho\uparrow1$.}
Since $q_n$ is even and nonincreasing on $[0,\sqrt n]$,
for every $|x|<\sqrt n$ we have
\[
 0\le V(\rho x)-V(0)\uparrow V(x)-V(0)
 \quad\text{as }\rho\uparrow1.
\]
For every coefficient vector with
$\sum\limits_{i=1}^n a_i^2=1$, the Cauchy--Schwarz inequality gives
\[
 \left|\sum\limits_{i=1}^n a_iU_i\right|
 \le\left(\sum\limits_{i=1}^n U_i^2\right)^{\frac12}
 <\sqrt n
 \quad\text{almost surely}.
\]
Thus both random sums in
\eqref{eq:cube-contracted-comparison} lie in the interval
where this monotone limit holds.

The density $q_n$ is bounded and supported on
$[-\sqrt n,\sqrt n]$, so
\[
 \int_{-\sqrt n}^{\sqrt n}|q_n(x)\log q_n(x)|\,\mathrm{d}x
 \le2\sqrt n\max\limits_{0\le u\le q_n(0)}|u\log u|<\infty.
\]
Subtracting $V(0)$ from both sides of
\eqref{eq:cube-contracted-comparison} and applying the monotone
convergence theorem now gives
\begin{equation}\label{eq:cube-cross-entropy}
 \mathbb E\left[V\left(\sum\limits_{i=1}^n a_iU_i\right)\right]
 \le\mathbb E\left[V\left(\frac1{\sqrt n}\sum\limits_{i=1}^n U_i\right)\right]
 =h\left(\frac1{\sqrt n}\sum\limits_{i=1}^n U_i\right)<\infty.
\end{equation}

\medskip
\noindent\emph{Entropy and equality.}
To obtain the entropy bound from \eqref{eq:cube-cross-entropy},
it remains to prove
\[
 h\left(\sum\limits_{i=1}^n a_iU_i\right)
 \le\mathbb E\left[V\left(\sum\limits_{i=1}^n a_iU_i\right)\right].
\]
We first verify that the entropy is finite, then establish
this comparison and determine when equality holds.

Since $\sum\limits_{i=1}^n a_i^2=1$, at least one coefficient is at
least $1/\sqrt n$. Relabeling if necessary, assume that
$a_1\ge1/\sqrt n$. Since $U_1$ is independent of
$U_2,\ldots,U_n$, conditional on these variables the sum
$\sum\limits_{i=1}^n a_iU_i$ is uniform on
\[
 \left[\sum\limits_{i=2}^n a_iU_i-a_1,\,
       \sum\limits_{i=2}^n a_iU_i+a_1\right].
\]
For every real $x$, its conditional density is therefore
\[
 \frac1{2a_1}
 \mathbf1_{\left\{\left|x-\sum\limits_{i=2}^n a_iU_i\right|\le a_1\right\}}.
\]
Taking the expectation of this conditional density gives a
density $f$ of $\sum\limits_{i=1}^n a_iU_i$, with
\[
 f(x)=\frac1{2a_1}\mathbb E\left[
   \mathbf1_{\left\{\left|x-\sum\limits_{i=2}^n a_iU_i\right|\le a_1\right\}}\right]
 =\frac1{2a_1}\mathbb P\left(
   \left|x-\sum\limits_{i=2}^n a_iU_i\right|\le a_1\right)
 \le\frac1{2a_1}\le\frac{\sqrt n}{2}.
\]
The inequalities $|U_i|\le1$, the nonnegativity of the coefficients, and the
Cauchy--Schwarz inequality give, almost surely,
\[
 \left|\sum\limits_{i=1}^n a_iU_i\right|
 \le\sum\limits_{i=1}^n a_i|U_i|
 \le\sum\limits_{i=1}^n a_i
 \le\sqrt{n\sum\limits_{i=1}^n a_i^2}=\sqrt n.
\]
Hence $f(x)=0$ for $|x|>\sqrt n$. The bound
$0\le f(x)\le\sqrt n/2$ also implies, for every real $x$,
\[
 |f(x)\log f(x)|
 \le\max\limits_{0\le u\le\frac{\sqrt n}{2}}|u\log u|<\infty.
\]
Integrating over
$[-\sqrt n,\sqrt n]$ therefore gives
\begin{align*}
 \int_{-\infty}^{\infty}|f(x)\log f(x)|\,\mathrm{d}x
 &=\int_{-\sqrt n}^{\sqrt n}|f(x)\log f(x)|\,\mathrm{d}x
 \le2\sqrt n\max\limits_{0\le u\le\frac{\sqrt n}{2}}|u\log u|<\infty.
\end{align*}
By \eqref{eq:cube-cross-entropy} and $V\ge V(0)$, we also have
\[
 \int_{-\sqrt n}^{\sqrt n}f(x)|\log q_n(x)|\,\mathrm{d}x
 =\mathbb E\left[\left|V\left(\sum\limits_{i=1}^n a_iU_i\right)\right|\right]
 \le\mathbb E\left[V\left(\sum\limits_{i=1}^n a_iU_i\right)-V(0)\right]+|V(0)|<\infty.
\]

We now compare the entropy with the expectation in
\eqref{eq:cube-cross-entropy}.
For $u\ge0$, the inequality $u\log u-u+1\ge0$ is strict
unless $u=1$. Since $q_n$ is positive on
$(-\sqrt n,\sqrt n)$, applying this inequality with
$u=f(x)/q_n(x)$ gives, for $|x|<\sqrt n$,
$f(x)\log[f(x)/q_n(x)]-f(x)+q_n(x)\ge0$.
Since $f$ and $q_n$ are probability densities supported in
$[-\sqrt n,\sqrt n]$, we have
\[
 \int_{-\sqrt n}^{\sqrt n}[-f(x)+q_n(x)]\,\mathrm{d}x
 =-\int_{-\sqrt n}^{\sqrt n}f(x)\,\mathrm{d}x
   +\int_{-\sqrt n}^{\sqrt n}q_n(x)\,\mathrm{d}x
 =-1+1=0.
\]
Since $f\log f$ and $f\log q_n$ are integrable, integrating
the pointwise inequality gives
\begin{equation}\label{eq:cube-relative-entropy}
 0\le\int_{-\sqrt n}^{\sqrt n}f(x)\log\frac{f(x)}{q_n(x)}\,\mathrm{d}x
 =-h\left(\sum\limits_{i=1}^n a_iU_i\right)
   +\mathbb E\left[V\left(\sum\limits_{i=1}^n a_iU_i\right)\right].
\end{equation}
Combining \eqref{eq:cube-relative-entropy} with
\eqref{eq:cube-cross-entropy} gives
\[
 h\left(\sum\limits_{i=1}^n a_iU_i\right)
 \le\mathbb E\left[V\left(\sum\limits_{i=1}^n a_iU_i\right)\right]
 \le h\left(\frac1{\sqrt n}\sum\limits_{i=1}^n U_i\right).
\]

If equality holds in the entropy bound, equality must also
hold in both \eqref{eq:cube-cross-entropy} and
\eqref{eq:cube-relative-entropy}. In particular, the integral
in \eqref{eq:cube-relative-entropy} is zero. Since both densities
have integral one, the nonnegative function $f\log(f/q_n)-f+q_n$
also has integral zero and therefore vanishes almost everywhere.
Since $u\log u-u+1=0$ only for $u=1$, this gives
$f=q_n$ almost everywhere.

Since $|U_i|\le1$ and $a_i\ge0$, we have $f(x)=0$
for $|x|>\sum\limits_{i=1}^n a_i$.
If $\sum\limits_{i=1}^n a_i<\sqrt n$, positivity of $q_n$ on
$(\sum\limits_{i=1}^n a_i,\sqrt n)$ gives the contradiction
\[
 0=\int_{\sum\limits_{i=1}^n a_i}^{\sqrt n}f(x)\,\mathrm{d}x
 =\int_{\sum\limits_{i=1}^n a_i}^{\sqrt n}q_n(x)\,\mathrm{d}x>0.
\]
The Cauchy--Schwarz inequality also gives
$\sum\limits_{i=1}^n a_i\le\sqrt n$, so
$\sum\limits_{i=1}^n a_i=\sqrt n$.
Equality in the Cauchy--Schwarz inequality now implies $a_i=1/\sqrt n$ for
every $1\le i\le n$.
Conversely, if $|a_i|=1/\sqrt n$ for every $1\le i\le n$, symmetry of
the $U_i$ makes the two sums in the proposition identically
distributed, and equality follows.
\end{proof}

\section{Exponential tilting}
\label{sec:q3-tilt-reduction}

Fix an integer $n\ge150$. By Proposition~\ref{prop:cube-criterion},
we need to prove the derivative inequality \eqref{eq:cube-third}.
For $1\le j\le n$, let $p_j$ be the density of
$\sum\limits_{i=1}^j U_i$. Since
$q_n(x)=\sqrt n\,p_n(\sqrt n\,x)$, taking logarithms and
differentiating three times gives, for $0<x<\sqrt n$,
\begin{equation}\label{eq:q3-density-rescaling}
 (\log q_n)'''(x)=n^{\frac32}(\log p_n)'''(\sqrt n\,x).
\end{equation}
Thus it suffices to prove that, for every $0<x<n$, $(\log p_n)'''(x)<0$.

Let $U$ be uniform on $[-1,1]$. For $t\in\mathbb R$, define
\begin{equation}\label{eq:q3-log-mgf}
 K(t)=\log\mathbb E[\exp(tU)]
 =\log\left[\frac12\int_{-1}^1 \exp(tu)\,\mathrm{d}u\right]
 =\log\frac{\sinh t}{t},
\end{equation}
where the value at zero is defined by continuity.
To obtain the \emph{exponential tilt} with parameter $t$, we multiply
the uniform density by $\exp(tu)$ and divide by $\mathbb E[\exp(tU)]$.
Let $U_t$ be a random variable with density
\[
 u\longmapsto
 \frac12\frac{\exp(tu)}{\mathbb E[\exp(tU)]}\mathbf1_{[-1,1]}(u)
 =\frac12 \exp\left[tu-K(t)\right]\mathbf1_{[-1,1]}(u).
\]
The parameter $t=0$ gives the original uniform density.
This family is also discussed in \cite{Keady2015,Quijano2018}.
For every real $s$, integration against this density gives
\begin{equation}\label{eq:q3-tilted-mgf}
 \mathbb E[\exp(sU_t)]
 =\frac12 \exp\left[-K(t)\right]\int_{-1}^1\exp\left[(t+s)u\right]\,\mathrm{d}u
 =\exp[K(t+s)-K(t)].
\end{equation}
Differentiating \eqref{eq:q3-tilted-mgf} at $s=0$ gives
$\mathbb E[U_t]=K'(t)$ and $\operatorname{Var}(U_t)=K''(t)>0$.
We center $U_t$ and normalize its variance by letting
\[
 Y=\frac{U_t-\mathbb E[U_t]}{\sqrt{\operatorname{Var}(U_t)}}.
\]
For every real $s$, equation~\eqref{eq:q3-tilted-mgf} gives
\begin{align*}
 \mathbb E[\exp(sY)]
 &=\exp\left[K\left(t+\frac{s}{\sqrt{\operatorname{Var}(U_t)}}\right)
             -K(t)-\frac{s\mathbb E[U_t]}{\sqrt{\operatorname{Var}(U_t)}}\right].
\end{align*}
In particular, $\mathbb E[Y]=0$, $\mathbb E[Y^2]=1$, and
$\mathbb E[Y^3]=K'''(t)/[\operatorname{Var}(U_t)]^{3/2}$.

Fix $t>0$, let $Y_1,\ldots,Y_n$ be independent copies of $Y$, and let $w$
be the density of
\[
 W=\frac1{\sqrt n}\sum\limits_{i=1}^n Y_i.
\]
The dependence of $Y$, $W$, and $w$ on $t$ is left implicit.
To determine the sign of
$(\log p_n)'''(n\mathbb E[U_t])$, we will express it in terms
of $(\log w)'''(0)$. We first compare the densities, using
convolution to add one summand at a time.
For $2\le j\le n$, conditional on
$\sum\limits_{i=1}^{j-1}U_i=u$, the sum
$\sum\limits_{i=1}^j U_i$ is uniform on $[u-1,u+1]$.
Its conditional density at $x$ is
$\frac12\mathbf1_{\{\lvert x-u\rvert\le1\}}$. Averaging with respect to
$p_{j-1}$ gives
\begin{equation}\label{eq:q3-uniform-convolution}
 p_j(x)
 =\frac12\int_{-\infty}^{\infty}
 p_{j-1}(u)\mathbf1_{\{\lvert x-u\rvert\le1\}}\,\mathrm du
 =\frac12\int_{x-1}^{x+1}p_{j-1}(u)\,\mathrm du.
\end{equation}

For the tilted sums, we claim that the density of a sum of
$j$ independent copies of $U_t$ at $x$ is
$\exp[tx-jK(t)]p_j(x)$.
We prove the claim by induction. For $j=1$, it follows from
the definition of $U_t$,
since $p_1(x)=\frac12\mathbf1_{[-1,1]}(x)$.
For $2\le j\le n$, suppose the sum of the first $j-1$ copies
has density $\exp[tu-(j-1)K(t)]p_{j-1}(u)$ at $u$.
Conditional on their sum being $u$, the sum of all $j$ copies
has density
$\frac12\exp[t(x-u)-K(t)]\mathbf1_{\{\lvert x-u\rvert\le1\}}$ at $x$.
Averaging over $u$ gives
\begin{align*}
 &\frac12\int_{x-1}^{x+1}
 \exp[tu-(j-1)K(t)]p_{j-1}(u)
 \exp[t(x-u)-K(t)]\,\mathrm du\\
 ={}&\exp[tx-jK(t)]\frac12\int_{x-1}^{x+1}p_{j-1}(u)\,\mathrm du\\
 ={}&\exp[tx-jK(t)]p_j(x),
\end{align*}
where the last equality follows from
\eqref{eq:q3-uniform-convolution}.
Taking $j=n$, the density of
the sum of $n$ independent copies of $U_t$ at $x$ is
$\exp[tx-nK(t)]p_n(x)$.

We now express this density through $w$.
The variables $\mathbb E[U_t]+\sqrt{\operatorname{Var}(U_t)}Y_i$
are independent copies of $U_t$, and
\[
 \sum\limits_{i=1}^n
 \left[\mathbb E[U_t]+\sqrt{\operatorname{Var}(U_t)}Y_i\right]
 =n\mathbb E[U_t]+\sqrt{n\operatorname{Var}(U_t)}W.
\]
Since $W$ has density $w$, the density of
$n\mathbb E[U_t]+\sqrt{n\operatorname{Var}(U_t)}W$ at $x$ is
\begin{equation}\label{eq:q3-affine-sum-density}
 \frac1{\sqrt{n\operatorname{Var}(U_t)}}
 w\left(\frac{x-n\mathbb E[U_t]}{\sqrt{n\operatorname{Var}(U_t)}}\right).
\end{equation}
Comparing \eqref{eq:q3-affine-sum-density} with
$\exp[tx-nK(t)]p_n(x)$, for $-n<x<n$ we obtain
\begin{equation}\label{eq:q3-tilt-density}
 p_n(x)
 =\frac{\exp\left[nK(t)-tx\right]}{\sqrt{n\operatorname{Var}(U_t)}}
   w\left(\frac{x-n\mathbb E[U_t]}{\sqrt{n\operatorname{Var}(U_t)}}\right).
\end{equation}

Before taking logarithms in \eqref{eq:q3-tilt-density}, we check
positivity and smoothness using \eqref{eq:q3-uniform-convolution}.
The density $p_2$ is continuous and positive on $(-2,2)$.
For $3\le j\le n$, differentiating
\eqref{eq:q3-uniform-convolution} gives
\[
 p_j'(x)=\frac12\left[p_{j-1}(x+1)-p_{j-1}(x-1)\right].
\]
Starting from $p_2$, each convolution adds one continuous
derivative, so $p_n$ has $n-2$ continuous derivatives.
For $|x|<j$, the intervals $(x-1,x+1)$ and $(-(j-1),j-1)$
overlap in an interval of positive length. If $p_{j-1}>0$
on $(-(j-1),j-1)$, then \eqref{eq:q3-uniform-convolution}
gives $p_j(x)>0$. Induction from $p_2$ therefore gives
$p_n>0$ on $(-n,n)$.

Since $n\ge150$, the regularity of $p_n$ and
\eqref{eq:q3-tilt-density} show that $w$ is positive and has
three continuous derivatives in the interior of its support.
Taking logarithms in \eqref{eq:q3-tilt-density}, for $-n<x<n$
we obtain
\[
 \log p_n(x)
 =nK(t)-tx-\frac12\log[n\operatorname{Var}(U_t)]
 +\log w\left(
       \frac{x-n\mathbb E[U_t]}{\sqrt{n\operatorname{Var}(U_t)}}\right).
\]
Differentiating three times with $t$ fixed and evaluating at
$x=n\mathbb E[U_t]$ gives
\begin{equation}\label{eq:q3-tilt-log-third}
 (\log p_n)'''\left(n\mathbb E[U_t]\right)
 =\frac{(\log w)'''(0)}{[n\operatorname{Var}(U_t)]^{\frac32}}.
\end{equation}
We derived \eqref{eq:q3-tilt-log-third} with $t$ fixed. We now
vary $t$ to cover all points of $(0,n)$.
The mean of $U_t$ is strictly increasing in $t$ because
$\frac{\mathrm{d}}{\mathrm{d}t}\mathbb E[U_t]
=K''(t)=\operatorname{Var}(U_t)>0$.
For $t>0$,
\begin{equation}\label{eq:q3-tilted-mean-range}
 \mathbb E[U_t]=\coth t-\frac1t,\qquad
 \lim_{t\downarrow0}\mathbb E[U_t]=0,\qquad
 \lim_{t\to\infty}\mathbb E[U_t]=1.
\end{equation}
Thus $t\mapsto n\mathbb E[U_t]$ maps $(0,\infty)$ onto $(0,n)$,
and for every $x\in(0,n)$ there is a unique $t>0$ with
$n\mathbb E[U_t]=x$. This choice of parameter is used in
saddlepoint approximation
\cite{Daniels1954}.
Equations
\eqref{eq:q3-density-rescaling} and \eqref{eq:q3-tilt-log-third}
reduce \eqref{eq:cube-third} to proving
$(\log w)'''(0)<0$ for every $t>0$.

To prove this sign, we separate the third density derivative
from the lower derivatives. Differentiating
$(\log w)''=w''/w-(w'/w)^2$ gives
\begin{equation}\label{eq:q3-tilt-jetnumerator}
 (\log w)'''(0)
 =\frac{w'''(0)}{w(0)}+\frac{w'(0)}{w(0)}
 \left[-\frac{3w''(0)}{w(0)}+\frac{2w'(0)^2}{w(0)^2}\right].
\end{equation}
For the standard Gaussian density, the expression in square
brackets in \eqref{eq:q3-tilt-jetnumerator} equals $3$.
In a formal Gaussian approximation, the third-moment correction
is proportional to $(x^3-3x)\exp(-x^2/2)$, whose first and third
derivatives at zero are $-3$ and $15$.
This suggests $w'''(0)\approx-5w'(0)$, and hence a comparison
constant strictly between $3$ and $5$. We choose $7/2$.

Thus it remains to establish the following three estimates
for every $t>0$:
\[
 w'(0)>0,\qquad
 -\frac{3w''(0)}{w(0)}+\frac{2w'(0)^2}{w(0)^2}<\frac72,
 \qquad w'''(0)+\frac72w'(0)<0.
\]
Since $w(0)>0$ and $w'(0)>0$, substitution in
\eqref{eq:q3-tilt-jetnumerator} gives
\[
 (\log w)'''(0)
 <\frac{w'''(0)+\frac72w'(0)}{w(0)}<0.
\]
The first two estimates are proved in
Proposition~\ref{prop:q3-even-jets}, and the third is proved in
Proposition~\ref{prop:q3-tilt-direct}.

\section{Sine and cosine estimates}
\label{sec:q3-global-odd-tail}

Fix an integer $n\ge150$ and $t>0$, and let $U_t$, $Y$, $W$, and $w$
be as in Section~\ref{sec:q3-tilt-reduction}. We first record bounds on
the variance of $U_t$ and the odd moments of $Y$. Their proof is given in
Appendix~\ref{sec:q3-even-moments80}.

\begin{lemma}\label{lem:q3-moments}
For every $t>0$,
\begin{equation}\label{eq:q3odd:moment-bounds}
 0<\operatorname{Var}(U_t)\le\frac13,\qquad -2\le\mathbb E[Y^3]<0,\qquad
 0<\frac{\mathbb E[Y^5]}{\mathbb E[Y^3]}<22.
\end{equation}
\end{lemma}

\subsection{Integral formulas for the density derivatives}

Our goal is to express $w(0)$, $w'(0)$, $w''(0)$, and $w'''(0)$ as
integrals involving $\mathbb E[\cos(uW)]$ and $\mathbb E[\sin(uW)]$.
Near $u=0$ we compare the integrands with Gaussian terms. To justify
inversion and bound the tails, we begin with an estimate for one summand.

\begin{lemma}\label{lem:q3-tilt-gamma}
For every $t>0$ and every real $v$,
\begin{equation}\label{eq:q3-real-transform-bound}
 \mathbb E[\cos(vY)]^2+\mathbb E[\sin(vY)]^2
 \le\frac1{1+v^2}.
\end{equation}
\end{lemma}

\begin{proof}
We first establish the scalar inequality used to bound the two expectations.
For $|s|\le\pi$, the bound in \cite[proof of Theorem~2]{AbiKhuzam2017} gives
\[
 \sin^2s\le s^2\exp\left(-\frac{s^2}{3}\right)\le\frac{s^2}{1+\frac13s^2}.
\]
For $|s|>\pi$, the last fraction exceeds one and hence still bounds $\sin^2s$.
Since $\operatorname{Var}(U_t)\le1/3$ by Lemma~\ref{lem:q3-moments}, we have
$[1+s^2\operatorname{Var}(U_t)]\sin^2s\le s^2$ for every real $s$.
Differentiating \eqref{eq:q3-log-mgf} twice gives
$\operatorname{Var}(U_t)=K''(t)=t^{-2}-(\sinh t)^{-2}$. Thus
\begin{equation}\label{eq:q3-sine-rational}
 \bigl[1+s^2\operatorname{Var}(U_t)\bigr]
 \left(1+\frac{\sin^2s}{\sinh^2t}\right)
 \le 1+s^2\left(\operatorname{Var}(U_t)+\frac{1}{\sinh^2t}\right)
 =1+\frac{s^2}{t^2}.
\end{equation}

We now compute the expectations in \eqref{eq:q3-real-transform-bound}.
Fix $v\in\mathbb R$ and let $s=v/\sqrt{\operatorname{Var}(U_t)}$.
Then $vY=s(U_t-\mathbb E[U_t])$ and $s^2\operatorname{Var}(U_t)=v^2$.
By \eqref{eq:q3-log-mgf}, $U_t$ has density $t\exp(tu)/(2\sinh t)$ on $[-1,1]$, so
\begin{align}
 \mathbb E[\cos(sU_t)]
 &=\frac{t}{2\sinh t}
       \int_{-1}^1\exp(tu)\cos(su)\,\mathrm du\notag\\
 &=\frac{t}{2(t^2+s^2)\sinh t}
       \left[\exp(tu)\bigl(t\cos(su)+s\sin(su)\bigr)\right]_{u=-1}^{u=1}\notag\\
 &=\frac{t}{t^2+s^2}(t\cos s+s\coth t\sin s),
       \label{eq:q3-tilted-cosine}\\[4pt]
 \mathbb E[\sin(sU_t)]
 &=\frac{t}{2\sinh t}
       \int_{-1}^1\exp(tu)\sin(su)\,\mathrm du\notag\\
 &=\frac{t}{2(t^2+s^2)\sinh t}
       \left[\exp(tu)\bigl(t\sin(su)-s\cos(su)\bigr)\right]_{u=-1}^{u=1}\notag\\
 &=\frac{t}{t^2+s^2}(t\coth t\sin s-s\cos s).
       \label{eq:q3-tilted-sine}
\end{align}
The trigonometric identities give the first equality below, and the mixed
terms cancel in the second. Substituting \eqref{eq:q3-tilted-cosine} and
\eqref{eq:q3-tilted-sine}, then applying \eqref{eq:q3-sine-rational}, gives
\begin{align}
 \mathbb E[\cos(vY)]^2
  +\mathbb E[\sin(vY)]^2
 &=\Bigl[\cos(s\mathbb E[U_t])\mathbb E[\cos(sU_t)]
          +\sin(s\mathbb E[U_t])\mathbb E[\sin(sU_t)]\Bigr]^2\notag\\
 &\phantom{=}\quad+\Bigl[\cos(s\mathbb E[U_t])\mathbb E[\sin(sU_t)]
          -\sin(s\mathbb E[U_t])\mathbb E[\cos(sU_t)]\Bigr]^2\notag\\
 &=\mathbb E[\cos(sU_t)]^2+\mathbb E[\sin(sU_t)]^2\notag\\
 &=\frac{t^2}{t^2+s^2}\left(1+\frac{\sin^2s}{\sinh^2t}\right)\label{eq:q3-tilt-modulus}\\
 &\le\frac{t^2}{t^2+s^2}\,
       \frac{t^2+s^2}{t^2[1+s^2\operatorname{Var}(U_t)]}\notag\\
 &=\frac{1}{1+s^2\operatorname{Var}(U_t)}\notag\\
 &=\frac{1}{1+v^2}.\notag
\end{align}
\end{proof}

We next pass from $Y$ to the sum $W$. For every real $u$, let
\begin{equation}\label{eq:q3-even-comparison}
 \rho(u)=\sqrt{\mathbb E[\cos(uW)]^2+\mathbb E[\sin(uW)]^2}.
\end{equation}
In particular, $\rho(u)$ bounds the absolute value of each expectation.
For independent real random variables $X$ and $Z$, the trigonometric
identities give
\begin{align}
 &\mathbb E[\cos(X+Z)]^2+\mathbb E[\sin(X+Z)]^2\notag\\
 ={}&\Bigl(\mathbb E[\cos X]\mathbb E[\cos Z]-\mathbb E[\sin X]\mathbb E[\sin Z]\Bigr)^2
    +\Bigl(\mathbb E[\sin X]\mathbb E[\cos Z]+\mathbb E[\cos X]\mathbb E[\sin Z]\Bigr)^2\notag\\
 ={}&\Bigl(\mathbb E[\cos X]^2+\mathbb E[\sin X]^2\Bigr)
     \Bigl(\mathbb E[\cos Z]^2+\mathbb E[\sin Z]^2\Bigr).
 \label{eq:q3-independent-trig}
\end{align}
For a single summand, \eqref{eq:q3-real-transform-bound} at $v=u/\sqrt n$
bounds the sum of squares by $(1+u^2/n)^{-1}$.
Since $uW=(u/\sqrt n)\sum\limits_{i=1}^nY_i$ and the $Y_i$ are independent
copies of $Y$, iterating \eqref{eq:q3-independent-trig} gives a product of
$n$ identical factors. We will use the same comparison function for the
cosine and sine estimates. For every real $u$, let
\begin{align*}
 r_n(u)&=\left(1+\frac{u^2}{n}\right)^{-\frac{n-1}{2}}.
\end{align*}
Taking the square root of the product gives
\begin{equation}
 \rho(u)
 =\left(\mathbb E\left[\cos\left(\frac{uY}{\sqrt n}\right)\right]^2
          +\mathbb E\left[\sin\left(\frac{uY}{\sqrt n}\right)\right]^2\right)^{\frac n2}
 \le\left(1+\frac{u^2}{n}\right)^{-\frac n2}\le r_n(u).
 \label{eq:q3-length-bound}
\end{equation}
For each $j\in\{0,1,2,3\}$, \eqref{eq:q3-length-bound} gives
\[
 \int_0^\infty u^j\rho(u)\,\mathrm{d}u
 \le\int_0^\infty u^jr_n(u)\,\mathrm{d}u<\infty.
\]
The last integral is finite because $r_n$ is bounded near zero
and $u^jr_n(u)=O(u^{j-n+1})$ as $u\to\infty$.
By \eqref{eq:q3-even-comparison}, the characteristic function of $W$ has
even modulus $\rho$ and is therefore integrable on $\mathbb R$.
For every $x\in\mathbb R$, we apply Fourier inversion to the density $w$
of $W$. Taking real parts and using that $\cos(u(W-x))$ is even in $u$ gives
\begin{align*}
 w(x)
 &=\frac{1}{2\pi}\int_{-\infty}^{\infty}\exp(-iux)\,\mathbb E[\exp(iuW)]\,\mathrm{d}u\\
 &=\frac{1}{2\pi}\int_{-\infty}^{\infty}\mathbb E[\exp(iu(W-x))]\,\mathrm{d}u\\
 &=\frac{1}{2\pi}\int_{-\infty}^{\infty}\mathbb E[\cos(u(W-x))]\,\mathrm{d}u\\
 &=\frac{1}{\pi}\int_0^\infty\cos(ux)\mathbb E[\cos(uW)]
 +\sin(ux)\mathbb E[\sin(uW)]\,\mathrm{d}u.
\end{align*}
For every real $x$, $u\ge0$, and $j\in\{0,1,2,3\}$, the
Cauchy--Schwarz inequality in $\mathbb R^2$ gives
\[
 \left|\frac{\mathrm{d}^j}{\mathrm{d}x^j}
 [\cos(ux)\mathbb E[\cos(uW)]+\sin(ux)\mathbb E[\sin(uW)]]\right|
 \le u^j\rho(u).
\]
Since $u^j\rho(u)$ is integrable and does not depend on $x$,
the dominated convergence theorem permits three
differentiations under the integral. Evaluating at $x=0$ gives
\begin{align}
 w(0)&=\frac1\pi\int_0^\infty\mathbb E[\cos(uW)]\,\mathrm{d}u,\notag\\
 w'(0)&=\frac1\pi\int_0^\infty u\mathbb E[\sin(uW)]\,\mathrm{d}u,\notag\\
 w''(0)&=-\frac1\pi\int_0^\infty u^2\mathbb E[\cos(uW)]\,\mathrm{d}u,\notag\\
 w'''(0)&=-\frac1\pi\int_0^\infty u^3\mathbb E[\sin(uW)]\,\mathrm{d}u.
 \label{eq:q3-fourier-jets}
\end{align}
\subsection{Bounds for the sine expectation}

We first prove absolute-value bounds for the sine expectations of
$Y$ and $W$ for all real arguments, then a lower bound near zero.
In \eqref{eq:q3-fourier-jets}, both $w'(0)$ and $w'''(0)$
are expressed through $\mathbb E[\sin(uW)]$.
At $t=0$, the variables $Y$ and $W$ are symmetric, so
$\mathbb E[Y^3]=w'(0)=w'''(0)=0$.
To determine the signs for arbitrarily small $t>0$, we need
to keep the factor $|\mathbb E[Y^3]|/\sqrt n$ in our estimates
for the sine expectation and its remainders.

For the centered tilted density $f$, the difference $f(y)-f(-y)$
changes sign once on the positive half-line. The cubic bound below
follows from the density comparison of Denuit, Lef\`evre, and Shaked
\cite[Theorem~4.4]{Denuit1998} applied to $Y$ and $-Y$, together with
the moment formula and equality criterion of Mattner and Shevtsova
\cite[Theorem~4.2(a)(b)]{MattnerShevtsova2019}. We use the same
integral comparison to prove this bound and to estimate the fifth-order
remainder in Lemma~\ref{lem:q3-sine-lower80}.

\begin{lemma}\label{lem:q3-onecrossing-odd}
Let $\xi$ have density $f$, with $\mathbb E[\xi]=0$ and
$\mathbb E[|\xi|^3]<\infty$. Suppose that, for some $b\ge0$,
$f(y)-f(-y)$ is nonnegative for almost every $0<y<b$ and
nonpositive for almost every $y>b$.
Then $\mathbb E[\xi^3]\le0$, and for every real $v$,
\begin{equation}\label{eq:q3-global-cubic-imag}
 |\mathbb E[\sin(v\xi)]|\le\frac16|\mathbb E[\xi^3]|\,|v|^3.
\end{equation}
If $f(y)-f(-y)$ is nonzero on a set of positive measure and
$v\ne0$, the inequality is strict.
\end{lemma}

\begin{proof}
If $f$ is even almost everywhere, both expectations
$\mathbb E[\xi^3]$ and $\mathbb E[\sin(v\xi)]$ vanish.
Otherwise $b>0$, since for $b=0$ the integrand in
$\int_0^\infty y[f(y)-f(-y)]\,\mathrm{d}y=\mathbb E[\xi]=0$
would be nonpositive and strictly negative on a set of positive measure.

Let $h$ be an odd function on $\mathbb R$ whose restriction to $[0,\infty)$ is convex, with $\mathbb E[|h(\xi)|]<\infty$.
Since $\mathbb E[\xi]=0$, subtracting a multiple of $y$ does not
change the expectation. We subtract the line through $(0,0)$ and
$(b,h(b))$ so that the remaining function has the opposite sign
to $f(y)-f(-y)$ on each side of $b$. To see this, convexity and
$h(0)=0$ show that $h(y)/y$ is nondecreasing for $y>0$. Hence
$h(y)-yh(b)/b$ is nonpositive for $0<y<b$ and nonnegative
for $y>b$. By oddness of $h$,
\begin{equation}\label{eq:q3-odd-convex-comparison}
 \mathbb E[h(\xi)]
 =\int_0^\infty\left(h(y)-\frac{h(b)}{b}\,y\right)
      [f(y)-f(-y)]\,\mathrm{d}y\le0.
\end{equation}
If $h$ is strictly convex on the positive half-line, $h(y)/y$ is
strictly increasing, and the integral is strictly negative.

To bound the sine expectation in both directions by the cubic moment,
for $v>0$ take $h_\pm(y)=v^3y^3/6\pm\sin(vy)$.
Both functions are strictly convex for $y>0$, since
$h_\pm''(y)=v^2[vy\mp\sin(vy)]>0$. Applying
\eqref{eq:q3-odd-convex-comparison} with $h_+$ and $h_-$ gives
\[
 \frac{v^3}{6}\mathbb E[\xi^3]+|\mathbb E[\sin(v\xi)]|<0.
\]
This gives $\mathbb E[\xi^3]<0$ and the strict bound for $v>0$.
At $v=0$ both sides of \eqref{eq:q3-global-cubic-imag} vanish.
Oddness of the sine gives the bound for $v<0$.
\end{proof}

\begin{corollary}\label{cor:q3-tilt-global-odd}
For the tilted variable $Y$ of Section~\ref{sec:q3-tilt-reduction},
every $t>0$, every integer $n\ge1$, and every real $v$,
\begin{align*}
 &\left|\mathbb E\left[\sin\left(v\sum\limits_{i=1}^nY_i\right)\right]\right|
 \le\frac16n|\mathbb E[Y^3]|\,|v|^3(1+v^2)^{-\frac{n-1}{2}}.
\end{align*}
If $v\ne0$, the inequality is strict.
\end{corollary}

\begin{proof}
Let $f$ be the density of $Y$. Since
$Y=(U_t-\mathbb E[U_t])/\sqrt{\operatorname{Var}(U_t)}$ and
$U_t$ takes values in $[-1,1]$, $f$ vanishes outside
\begin{equation}\label{eq:q3-tilted-support}
 \left[-\frac{1+\mathbb E[U_t]}{\sqrt{\operatorname{Var}(U_t)}},
        \frac{1-\mathbb E[U_t]}{\sqrt{\operatorname{Var}(U_t)}}\right].
\end{equation}
On this interval, the change of variables
$u=\mathbb E[U_t]+\sqrt{\operatorname{Var}(U_t)}\,y$
in the density $\frac12\exp[tu-K(t)]$ of $U_t$ gives
\[
 f(y)=\frac12\sqrt{\operatorname{Var}(U_t)}
 \exp\left[t\left(\mathbb E[U_t]
       +\sqrt{\operatorname{Var}(U_t)}\,y\right)-K(t)\right].
\]
Since $t\sqrt{\operatorname{Var}(U_t)}>0$, the exponent is
strictly increasing in $y$, and hence so is $f$ on this interval.

Equation~\eqref{eq:q3-tilted-mean-range} gives $0<\mathbb E[U_t]<1$,
so the negative endpoint in~\eqref{eq:q3-tilted-support} has larger
absolute value than the positive endpoint. For positive $y$ below the
right endpoint, both $y$ and $-y$ lie in the interior of the interval, and strict monotonicity
gives $f(y)>f(-y)$. Above the right endpoint, $f(y)=0$.
We therefore obtain
\begin{align}
 0<y<\frac{1-\mathbb E[U_t]}{\sqrt{\operatorname{Var}(U_t)}}
 &\ \Longrightarrow\ f(y)-f(-y)>0,\notag\\
 y>\frac{1-\mathbb E[U_t]}{\sqrt{\operatorname{Var}(U_t)}}
 &\ \Longrightarrow\ f(y)-f(-y)\le0.
 \label{eq:q3-tilted-density-sign}
\end{align}
The variable $Y$ is bounded and centered, and the first inequality
in \eqref{eq:q3-tilted-density-sign} holds on an interval of positive
length. Lemma~\ref{lem:q3-onecrossing-odd} therefore proves the case
$n=1$, including strictness when $v\ne0$.
The case $v=0$ is immediate, so assume $v\ne0$.

For $n\ge2$, let
\[
 r=\sqrt{\mathbb E[\cos(vY)]^2+\mathbb E[\sin(vY)]^2}.
\]
Equation~\eqref{eq:q3-tilt-modulus} gives $r>0$, so we may choose
$\alpha\in[0,2\pi]$ such that
\[
 \mathbb E[\cos(vY)]=r\cos\alpha,\qquad
 \mathbb E[\sin(vY)]=r\sin\alpha.
\]
Applying the sine and cosine identities used in
\eqref{eq:q3-independent-trig} successively to the independent summands gives
$\mathbb E[\sin(v\sum\limits_{i=1}^nY_i)]=r^n\sin(n\alpha)$.
Using $|\sin(n\alpha)|\le n|\sin\alpha|$ and
$r\le(1+v^2)^{-1/2}$ from \eqref{eq:q3-real-transform-bound},
we obtain
\[
 \left|\mathbb E\left[\sin\left(v\sum\limits_{i=1}^nY_i\right)\right]\right|
 \le nr^n|\sin\alpha|
 =nr^{n-1}|\mathbb E[\sin(vY)]|
 \le n|\mathbb E[\sin(vY)]|(1+v^2)^{-\frac{n-1}{2}}.
\]
Lemma~\ref{lem:q3-onecrossing-odd} bounds the remaining sine
expectation and makes the result strict.
\end{proof}

We next compare the sine expectation with its cubic term.
The fifth moment will bound the error.

\begin{lemma}\label{lem:q3-sine-lower80}
For the tilted variable $Y$ of Section~\ref{sec:q3-tilt-reduction},
every $t>0$, and every $v\ge0$,
\[
 \mathbb E[\sin(vY)]
 \ge\frac16|\mathbb E[Y^3]|v^3\left(1-\frac{11v^2}{10}\right).
\]
\end{lemma}

\begin{proof}
We want to compare the sine expectation with the third and fifth
moments. To apply \eqref{eq:q3-odd-convex-comparison}, let
$\varphi(x)=x^5/120-x^3/6-\sin x$. Then $\varphi$ is convex on
$[0,\infty)$.
Thus $y\mapsto \varphi(vy)$ is odd and convex on the positive half-line
for $v\ge0$. The tilted variable $Y$ is bounded, and its density
satisfies \eqref{eq:q3-tilted-density-sign}. Hence
\eqref{eq:q3-odd-convex-comparison} gives $\mathbb E[\varphi(vY)]\le0$.
Since $\mathbb E[Y^3]<0$ and
$0<\mathbb E[Y^5]/\mathbb E[Y^3]<22$ by
\eqref{eq:q3odd:moment-bounds}, we obtain
\begin{align*}
 \mathbb E[\sin(vY)]
 &\ge-\frac{v^3}{6}\mathbb E[Y^3]+\frac{v^5}{120}\mathbb E[Y^5]
 \ge\frac16|\mathbb E[Y^3]|v^3\left(1-\frac{11v^2}{10}\right).
\end{align*}
\end{proof}

\subsection{Lower bounds for the cosine and sine expectations}

We next derive the lower bounds that will be integrated in the rest of
the proof. The cosine bound will estimate $w(0)$. The sine bound will
estimate $w'(0)$ and the signed integral in
Section~\ref{sec:q3-tilt-direct}.
By Lemma~\ref{lem:q3-moments}, $\mathbb E[Y^3]<0$, so for every real $u$ let
\begin{equation}\label{eq:q3-normalized-sine40}
 S_n(u)=\frac{\sqrt n}{|\mathbb E[Y^3]|}
 \mathbb E\left[\sin\left(\frac{u}{\sqrt n}\sum\limits_{i=1}^nY_i\right)\right]
 =\frac{\sqrt n}{|\mathbb E[Y^3]|}\mathbb E[\sin(uW)].
\end{equation}

\begin{lemma}\label{lem:q3-central-polar}
Let an integer $n\ge150$ and $t>0$ be given.
Then for every real $u$ with $|u|\le\sqrt n/2$,
\begin{equation}\label{eq:q3-aux80-central-cosine}
 \mathbb E[\cos(uW)]
 \ge\exp\left(-\frac12u^2\right)\left(1-\frac{u^4}{7n}\right)
       -\frac{u^6}{12n}r_n(u).
\end{equation}
For every $u$ with $0\le u\le\sqrt n/2$,
\begin{equation}\label{eq:q3-aux80-central-sine}
 S_n(u)
 \ge\frac{u^3}{6}\exp\left(-\frac12u^2\right)
       \left(1-\frac{11u^2}{10n}-\frac{u^4}{7n}\right)
       -\frac{u^9}{216n}r_n(u).
\end{equation}
\end{lemma}

\begin{proof}
Fix $0\le u\le\sqrt n/2$.
Since $\mathbb E[Y^2]=1$ and $u^2/n\le1/4$, we have
\begin{align}
 \mathbb E\left[\cos\left(\frac{uY}{\sqrt n}\right)\right]
 &\ge1-\frac{u^2}{2n}\ge\frac78>0.
 \label{eq:q3-positive-cosine}
\end{align}
Lemma~\ref{lem:q3-sine-lower80}, applied with parameter $v=u/\sqrt n$, gives
\[
 \mathbb E\left[\sin\left(\frac{uY}{\sqrt n}\right)\right]
 \ge\frac{|\mathbb E[Y^3]|u^3}{6n^{3/2}}
       \left(1-\frac{11u^2}{10n}\right)\ge0.
\]
We may therefore let
\[
 \theta(u)=n\arctan\left(
 \frac{\mathbb E\left[\sin\left(\frac{uY}{\sqrt n}\right)\right]}{\mathbb E\left[\cos\left(\frac{uY}{\sqrt n}\right)\right]}\right),
\]
with $0\le\theta(u)/n<\pi/2$.
By \eqref{eq:q3-length-bound}, we have
\[
 \mathbb E\left[\cos\left(\frac{uY}{\sqrt n}\right)\right]^2+\mathbb E\left[\sin\left(\frac{uY}{\sqrt n}\right)\right]^2
 =\rho(u)^{\frac2n}.
\]
Thus
\[
 \mathbb E\left[\cos\left(\frac{uY}{\sqrt n}\right)\right]=\rho(u)^{\frac1n}\cos\left(\frac{\theta(u)}n\right),
 \qquad
 \mathbb E\left[\sin\left(\frac{uY}{\sqrt n}\right)\right]=\rho(u)^{\frac1n}\sin\left(\frac{\theta(u)}n\right).
\]
Applying the trigonometric identities in
\eqref{eq:q3-independent-trig} to
$uW=\frac{u}{\sqrt n}\sum\limits_{i=1}^nY_i$ gives
\begin{equation}\label{eq:q3-real-polar}
 \mathbb E[\cos(uW)]=\rho(u)\cos\theta(u),\qquad
 \mathbb E[\sin(uW)]=\rho(u)\sin\theta(u).
\end{equation}

We first bound $\rho(u)$ and $\theta(u)$.
By \eqref{eq:q3-positive-cosine}, we have
\[
 \rho(u)\ge\mathbb E\left[\cos\left(\frac{uY}{\sqrt n}\right)\right]^n
 \ge\left(1-\frac{u^2}{2n}\right)^n.
\]
Since $0\le u^2/n\le1/4$, we have
\[
 n\log\left(1-\frac{u^2}{2n}\right)+\frac{u^2}{2}+\frac{u^4}{7n}
 =n\int_0^{\frac{u^2}{n}}\frac{s(1-4s)}{14(2-s)}\,\mathrm ds\ge0.
\]
Together with \eqref{eq:q3-length-bound}, this yields
\begin{equation}\label{eq:q3-direct-amplitude}
 \exp\left(-\frac12u^2-\frac{u^4}{7n}\right)
 \le\rho(u)\le r_n(u).
\end{equation}
For the angle, apply Lemma~\ref{lem:q3-onecrossing-odd} with
$\xi=Y$ and parameter $v=u/\sqrt n$.
Together with $\arctan x\le x$ for $x\ge0$ and
\eqref{eq:q3-positive-cosine}, this gives
\begin{align}
 \theta(u)
 &\le\frac{n\mathbb E\left[\sin\left(\frac{uY}{\sqrt n}\right)\right]}
           {\mathbb E\left[\cos\left(\frac{uY}{\sqrt n}\right)\right]}
 \le\frac{|\mathbb E[Y^3]|u^3}{6\sqrt n\left(1-\frac{u^2}{2n}\right)}
 \le\frac{|\mathbb E[Y^3]|u^3}{6\sqrt n\cdot\frac78}
 =\frac{4|\mathbb E[Y^3]|u^3}{21\sqrt n}.
 \label{eq:q3-even-phase}
\end{align}

We now prove the cosine bound. By \eqref{eq:q3odd:moment-bounds},
\eqref{eq:q3-direct-amplitude}, and \eqref{eq:q3-even-phase}, we have
\begin{align*}
 \frac12\rho(u)\theta(u)^2
 &\le\frac12r_n(u)\left(\frac{4|\mathbb E[Y^3]|u^3}{21\sqrt n}\right)^2
 \le\frac{32u^6}{441n}r_n(u)
 \le\frac{u^6}{12n}r_n(u).
\end{align*}
The inequality $\cos x\ge1-x^2/2$ and
\eqref{eq:q3-direct-amplitude} therefore give
\[
 \mathbb E[\cos(uW)]
 \ge\rho(u)-\frac12\rho(u)\theta(u)^2
 \ge\exp\left(-\frac12u^2-\frac{u^4}{7n}\right)
       -\frac{u^6}{12n}r_n(u).
\]
Applying $\exp(-x)\ge1-x$ to the factor
$\exp(-u^4/(7n))$ proves \eqref{eq:q3-aux80-central-cosine}.
Both sides of this inequality are even in $u$, so it also holds
for $-\sqrt n/2\le u<0$.

For the sine bound, since $\theta(u)\ge0$, the inequality
$\sin x\ge x-x^3/6$ for $x\ge0$ and
\eqref{eq:q3-real-polar} give
\begin{equation}\label{eq:q3-sine-cubic-decomposition}
 \mathbb E[\sin(uW)]
 \ge\rho(u)\theta(u)-\frac16\rho(u)\theta(u)^3.
\end{equation}
We now bound the first term from below and the second term from above.
For the first term, we have
\begin{align}
 \theta(u)
 &\ge n\rho(u)^{\frac1n}\sin\left(\frac{\theta(u)}n\right)
 =n\mathbb E\left[\sin\left(\frac{uY}{\sqrt n}\right)\right]
 \ge\frac{|\mathbb E[Y^3]|u^3}{6\sqrt n}
       \left(1-\frac{11u^2}{10n}\right)\ge0,
 \label{eq:q3-aux80-phase}
\end{align}
where the first inequality follows from $\sin x\le x$ for $x\ge0$
and $\rho(u)^{1/n}\le1$, and the second inequality follows from
Lemma~\ref{lem:q3-sine-lower80}, applied with parameter $v=u/\sqrt n$.

\begin{samepage}
Multiplying the nonnegative lower bounds in
\eqref{eq:q3-direct-amplitude} and \eqref{eq:q3-aux80-phase}
and applying $\exp(-x)\ge1-x$ to the factor
$\exp(-u^4/(7n))$, we obtain
\begin{align*}
 \rho(u)\theta(u)
 &\ge\frac{|\mathbb E[Y^3]|u^3}{6\sqrt n}
       \exp\left(-\frac12u^2\right)
       \exp\left(-\frac{u^4}{7n}\right)
       \left(1-\frac{11u^2}{10n}\right)\\
 &\ge\frac{|\mathbb E[Y^3]|u^3}{6\sqrt n}
       \exp\left(-\frac12u^2\right)
       \left(1-\frac{u^4}{7n}\right)
       \left(1-\frac{11u^2}{10n}\right)\\
 &\ge\frac{|\mathbb E[Y^3]|u^3}{6\sqrt n}
       \exp\left(-\frac12u^2\right)
       \left(1-\frac{11u^2}{10n}-\frac{u^4}{7n}\right).
\end{align*}
\end{samepage}

For the cubic term, \eqref{eq:q3odd:moment-bounds},
\eqref{eq:q3-direct-amplitude}, and \eqref{eq:q3-even-phase} give
\begin{align*}
 \frac16\rho(u)\theta(u)^3
 &\le\frac16r_n(u)\left(\frac{4|\mathbb E[Y^3]|u^3}{21\sqrt n}\right)^3
 \le\frac{128|\mathbb E[Y^3]|u^9}{27783n\sqrt n}r_n(u)
 \le\frac{|\mathbb E[Y^3]|u^9}{216n\sqrt n}r_n(u).
\end{align*}
Substituting the two estimates into
\eqref{eq:q3-sine-cubic-decomposition} and multiplying by
$\sqrt n/|\mathbb E[Y^3]|$ proves
\eqref{eq:q3-aux80-central-sine}.
\end{proof}

\subsection{The density and its first two derivatives}
\label{sec:q3-tilt-even}

To obtain lower bounds on $w(0)$ and $w'(0)$, we split the
corresponding integrals in \eqref{eq:q3-fourier-jets} at $u=\sqrt n/2$.
On $0\le u\le\sqrt n/2$, we use \eqref{eq:q3-aux80-central-cosine}
in the cosine integral and \eqref{eq:q3-aux80-central-sine} in the sine integral.
For $u>\sqrt n/2$, we bound the integrands in absolute value.

For every integer $k\ge0$ and every real $n>2k+2$,
define the even moment integrals by
\begin{equation}\label{eq:q3-comparison-moments}
 M_{2k}(n)=\frac1{\sqrt{2\pi}}\int_{-\infty}^{\infty}
 u^{2k}\left(1+\frac{u^2}{n}\right)^{-\frac{n-1}{2}}\,\mathrm du
 =\frac1{\sqrt{2\pi}}\int_{-\infty}^{\infty}u^{2k}r_n(u)\,\mathrm du.
\end{equation}
The condition $n>2k+2$ ensures that the integral is finite.
The moment and tail estimates are collected in
Appendix~\ref{sec:q3-moment-tails80}.

\begin{proposition}\label{prop:q3-even-jets}
For every integer $n\ge150$ and every $t>0$,
\begin{equation}\label{eq:q3-aux40-values}
 w(0)>\frac{19}{20\sqrt{2\pi}},\qquad
 -w''(0)<\frac{21}{20\sqrt{2\pi}},\qquad
 0<w'(0)<\frac{|\mathbb E[Y^3]|}{\sqrt{2\pi n}}.
\end{equation}
Moreover,
\begin{equation}\label{eq:q3-aux80-score}
 -\frac{3w''(0)}{w(0)}+\frac{2w'(0)^2}{w(0)^2}<\frac72.
\end{equation}
\end{proposition}

\begin{proof}
We first bound the derivatives in absolute value.
Applying Corollary~\ref{cor:q3-tilt-global-odd} with $v=u/\sqrt n$
and using the definition of $S_n$ in \eqref{eq:q3-normalized-sine40}, we obtain
\begin{equation}\label{eq:q3-global-envelope40}
 |S_n(u)|\le\frac{|u|^3}{6}\left(1+\frac{u^2}{n}\right)^{-\frac{n-1}{2}}
 =\frac{|u|^3}{6}r_n(u).
\end{equation}
By \eqref{eq:q3-even-comparison} and \eqref{eq:q3-length-bound}, we have
\begin{align}\label{eq:q3-cosine-envelope}
 |\mathbb E[\cos(uW)]|&\le\rho(u)\le r_n(u).
\end{align}
Taking absolute values for $w''(0)$ in \eqref{eq:q3-fourier-jets}
and using \eqref{eq:q3-cosine-envelope}, evenness,
\eqref{eq:q3-comparison-moments}, and the moment bound
\eqref{eq:q3-aux80-moment-bounds} in Appendix~\ref{sec:q3-moment-tails80}, we obtain
\[
 \sqrt{2\pi}\,|w''(0)|
 \le\frac2{\sqrt{2\pi}}\int_0^\infty u^2r_n(u)\,\mathrm du
 =M_2(n)<\frac{21}{20}.
\]
For $w'(0)$, we substitute
$\mathbb E[\sin(uW)]=|\mathbb E[Y^3]|S_n(u)/\sqrt n$
from \eqref{eq:q3-normalized-sine40} into \eqref{eq:q3-fourier-jets}.
Using \eqref{eq:q3-global-envelope40}, evenness,
\eqref{eq:q3-comparison-moments}, and the moment bound
\eqref{eq:q3-aux80-moment-bounds} in Appendix~\ref{sec:q3-moment-tails80}, we obtain
\[
 \frac{\sqrt{2\pi n}\,|w'(0)|}{|\mathbb E[Y^3]|}
 \le\frac2{\sqrt{2\pi}}\int_0^\infty u|S_n(u)|\,\mathrm du
 \le\frac2{6\sqrt{2\pi}}\int_0^\infty u^4r_n(u)\,\mathrm du
 =\frac16M_4(n)<1.
\]
The cosine integrand in \eqref{eq:q3-fourier-jets} is even in $u$.
Writing the integral over $\mathbb R$ and splitting it at
$|u|=\sqrt n/2$, we obtain
\begin{align*}
 &\sqrt{2\pi}\,w(0)\\
 ={}&\frac1{\sqrt{2\pi}}\left(
 \int_{-\frac{\sqrt n}{2}}^{\frac{\sqrt n}{2}}\mathbb E[\cos(uW)]\,\mathrm du
 +\int_{|u|>\frac{\sqrt n}{2}}\mathbb E[\cos(uW)]\,\mathrm du
 \right)\\
 \ge{}&\frac1{\sqrt{2\pi}}\left(\int_{-\frac{\sqrt n}{2}}^{\frac{\sqrt n}{2}}
 \left[\exp\left(-\frac12u^2\right)\left(1-\frac{u^4}{7n}\right)
 -\frac{u^6r_n(u)}{12n}\right]\,\mathrm du
 -\int_{|u|>\frac{\sqrt n}{2}}r_n(u)\,\mathrm du\right)\\
 \ge{}&\frac1{\sqrt{2\pi}}\left(\int_{-\infty}^{\infty}
 \left[\exp\left(-\frac12u^2\right)\left(1-\frac{u^4}{7n}\right)
 -\frac{u^6r_n(u)}{12n}\right]\,\mathrm du
 -\int_{|u|>\frac{\sqrt n}{2}}
 \left[\exp\left(-\frac12u^2\right)+r_n(u)\right]\,\mathrm du\right)\\
 ={}&1-\frac3{7n}-\frac{M_6(n)}{12n}-\frac1{\sqrt{2\pi}}
   \int_{|u|>\frac{\sqrt n}{2}}
     \left[\exp\left(-\frac12u^2\right)+r_n(u)\right]\,\mathrm du,
\end{align*}
where the first inequality follows from \eqref{eq:q3-aux80-central-cosine}
and \eqref{eq:q3-cosine-envelope}, and the last equality follows from
\eqref{eq:q3-comparison-moments}
and \eqref{eq:q3-gaussian-moments}.
By $M_6(n)<24$ in \eqref{eq:q3-aux80-moment-bounds} and the tail bound
\eqref{eq:q3-tail-zero}, for $n\ge150$ we have
\[
 \sqrt{2\pi}\,w(0)
 >1-\frac3{7n}-\frac{24}{12n}-\frac1n>\frac{19}{20}.
\]

Using the definition of $S_n$ and \eqref{eq:q3-fourier-jets}, we obtain
\begin{align*}
 &\frac{\sqrt{2\pi n}\,w'(0)}{|\mathbb E[Y^3]|}\\
 ={}&\sqrt{\frac2\pi}\int_0^\infty uS_n(u)\,\mathrm du\\
 \ge{}&\sqrt{\frac2\pi}\left(
 \int_0^{\frac{\sqrt n}{2}}uS_n(u)\,\mathrm du
 -\frac16\int_{\frac{\sqrt n}{2}}^\infty u^4r_n(u)\,\mathrm du\right)\\
 \ge{}&\frac16\sqrt{\frac2\pi}\left(
 \int_0^{\frac{\sqrt n}{2}}u^4
 \left[\exp\left(-\frac12u^2\right)
 \left(1-\frac{11u^2}{10n}-\frac{u^4}{7n}\right)
 -\frac{u^6}{36n}r_n(u)\right]\,\mathrm du
 -\int_{\frac{\sqrt n}{2}}^\infty u^4r_n(u)\,\mathrm du\right)\\
 \ge{}&\frac36-\frac{11\cdot15}{60n}-\frac{105}{42n}
 -\frac{M_{10}(n)}{216n}-\frac1{6\sqrt{2\pi}}
 \int_{|u|>\frac{\sqrt n}{2}}
 u^4\left[\exp\left(-\frac12u^2\right)+r_n(u)\right]\,\mathrm du\\
 ={}&\frac12-\frac{21}{4n}-\frac{M_{10}(n)}{216n}-\frac1{6\sqrt{2\pi}}
 \int_{|u|>\frac{\sqrt n}{2}}
 u^4\left[\exp\left(-\frac12u^2\right)+r_n(u)\right]\,\mathrm du,
\end{align*}
where the first inequality follows from \eqref{eq:q3-global-envelope40}
and the second inequality follows from \eqref{eq:q3-aux80-central-sine}.
The integral evaluations use evenness,
\eqref{eq:q3-comparison-moments}, and \eqref{eq:q3-gaussian-moments}.
By $M_{10}(n)<1500$ in \eqref{eq:q3-aux80-moment-bounds} and the tail bound
\eqref{eq:q3-tail-four}, for $n\ge150$ we have
\[
 \frac{\sqrt{2\pi n}\,w'(0)}{|\mathbb E[Y^3]|}
 >\frac12-\frac{21}{4n}-\frac{1500}{216n}-\frac1{2n}>\frac25.
\]
This proves \eqref{eq:q3-aux40-values}.

By \eqref{eq:q3-aux40-values}, $-w''(0)/w(0)<21/19$ and
$0<w'(0)/w(0)<20|\mathbb E[Y^3]|/(19\sqrt n)$.
Using $|\mathbb E[Y^3]|\le2$ from \eqref{eq:q3odd:moment-bounds}
and $n\ge150$, we obtain
\[
 -\frac{3w''(0)}{w(0)}+\frac{2w'(0)^2}{w(0)^2}
 <\frac{63}{19}+\frac8n\left(\frac{20}{19}\right)^2<\frac72.
\]
\end{proof}

\subsection{The third logarithmic derivative}
\label{sec:q3-tilt-direct}

Section~\ref{sec:q3-tilt-reduction} reduces the entropy comparison
to the sign of $(\log w)'''(0)$.
The bounds \eqref{eq:q3-aux40-values} and \eqref{eq:q3-aux80-score},
substituted in \eqref{eq:q3-tilt-jetnumerator}, show that it suffices to prove
$w'''(0)+(7/2)w'(0)<0$.
We use the two sine integrals in \eqref{eq:q3-fourier-jets}
to estimate this combination.

\begin{proposition}\label{prop:q3-tilt-direct}
Let $n\ge150$ be an integer. For every $t>0$, the standardized tilted density $w$ of
Section~\ref{sec:q3-tilt-reduction} satisfies
\begin{equation}\label{eq:q3-direct-signed}
 w'''(0)+\frac72w'(0)<0.
\end{equation}
\end{proposition}

\begin{proof}
Fix $n\ge150$ and $t>0$.
Using the definition of $S_n$ in \eqref{eq:q3-normalized-sine40},
the two sine integrals in \eqref{eq:q3-fourier-jets}, and evenness, we obtain
\begin{align}
 &-\frac{\sqrt{2\pi n}}{|\mathbb E[Y^3]|}
 \left[w'''(0)+\frac72w'(0)\right]\notag\\
 ={}&\sqrt{\frac2\pi}\int_0^\infty
 u\left(u^2-\frac72\right)S_n(u)\,\mathrm du
 \label{eq:q3-signed-fourier-integral}\\
 ={}&\sqrt{\frac2\pi}\int_0^\infty u\left(u^2-\frac72\right)
 \left[S_n(u)-\frac{u^3}{6}\exp\left(-\frac12u^2\right)\right]\,\mathrm du\notag\\
 &\qquad+\frac1{6\sqrt{2\pi}}\int_{-\infty}^{\infty}
 u^6\exp\left(-\frac12u^2\right)\,\mathrm du
 -\frac7{12\sqrt{2\pi}}\int_{-\infty}^{\infty}
 u^4\exp\left(-\frac12u^2\right)\,\mathrm du\notag\\
 ={}&\sqrt{\frac2\pi}\int_0^\infty u\left(u^2-\frac72\right)
 \left[S_n(u)-\frac{u^3}{6}\exp\left(-\frac12u^2\right)\right]\,\mathrm du+\frac34,
 \label{eq:q3-signed-gaussian200}
\end{align}
where the last equality follows from \eqref{eq:q3-gaussian-moments}.
We split the integral in \eqref{eq:q3-signed-gaussian200}
at $u=\sqrt{14}/2$ and $u=\sqrt n/2$.

For $0\le u<\sqrt{14}/2$, \eqref{eq:q3-global-envelope40}
and \eqref{eq:q3-envelope-near-zero} give
\[
 S_n(u)\le\frac{u^3}{6}r_n(u)
 \le\frac{21}{20}\cdot\frac{u^3}{6}\exp\left(-\frac12u^2\right).
\]
Multiplying by the nonpositive weight and using
\eqref{eq:q3-negative-weight-loss200}, we obtain
\begin{align*}
 &\sqrt{\frac2\pi}\int_0^{\frac{\sqrt{14}}{2}}
 u\left(u^2-\frac72\right)
 \left[S_n(u)-\frac{u^3}{6}\exp\left(-\frac12u^2\right)\right]\,\mathrm du\\
 \ge{}&-\frac1{20}\sqrt{\frac2\pi}\int_0^{\frac{\sqrt{14}}{2}}
 u\left(\frac72-u^2\right)\frac{u^3}{6}
 \exp\left(-\frac12u^2\right)\,\mathrm du>-\frac1{20}.
\end{align*}

For $\sqrt{14}/2\le u\le\sqrt n/2$, \eqref{eq:q3-aux80-central-sine} and
\eqref{eq:q3-integrated-central-error} give
\begin{align}
 &\sqrt{\frac2\pi}\int_{\frac{\sqrt{14}}{2}}^{\frac{\sqrt n}{2}}
 u\left(u^2-\frac72\right)
 \left[S_n(u)-\frac{u^3}{6}\exp\left(-\frac12u^2\right)\right]\,\mathrm du\notag\\
 \ge{}&-\frac1n\sqrt{\frac2\pi}\int_{\frac{\sqrt{14}}{2}}^{\frac{\sqrt n}{2}}
 u\left(u^2-\frac72\right)
 \left[\left(\frac{11u^5}{60}+\frac{u^7}{42}\right)
 \exp\left(-\frac12u^2\right)+\frac{u^9}{216}r_n(u)\right]\,\mathrm du\notag\\
 \ge{}&-\frac1n\sqrt{\frac2\pi}\int_{\frac{\sqrt{14}}{2}}^\infty
 u\left(u^2-\frac72\right)
 \left[\left(\frac{11u^5}{60}+\frac{u^7}{42}\right)
 \exp\left(-\frac12u^2\right)+\frac{u^9}{216}r_n(u)\right]\,\mathrm du
 >-\frac{80}{n}.
 \label{eq:q3-central-loss200}
\end{align}

For $u>\sqrt n/2$, we have
$0<u\left(u^2-\frac72\right)\le u^3$.
Using \eqref{eq:q3-global-envelope40}, evenness, and
\eqref{eq:q3-tail-six}, we obtain
\begin{align*}
 &\sqrt{\frac2\pi}\int_{\frac{\sqrt n}{2}}^\infty
 u\left(u^2-\frac72\right)
 \left[S_n(u)-\frac{u^3}{6}\exp\left(-\frac12u^2\right)\right]\,\mathrm du\\
 \ge{}&-\sqrt{\frac2\pi}\int_{\frac{\sqrt n}{2}}^\infty
 u\left(u^2-\frac72\right)
 \left[\frac{u^3}{6}\exp\left(-\frac12u^2\right)+|S_n(u)|\right]\,\mathrm du\\
 \ge{}&-\frac1{6\sqrt{2\pi}}\int_{|u|>\frac{\sqrt n}{2}}
 u^6\left[\exp\left(-\frac12u^2\right)+r_n(u)\right]\,\mathrm du
 >-\frac1{10}.
\end{align*}

Adding the estimates on the three intervals and using
\eqref{eq:q3-signed-fourier-integral} and \eqref{eq:q3-signed-gaussian200},
we obtain
\begin{align}
 \sqrt{\frac2\pi}\int_0^\infty
 u\left(u^2-\frac72\right)S_n(u)\,\mathrm du
 &>\frac34-\frac1{20}-\frac{80}{n}-\frac1{10}>0.
 \label{eq:q3-signed-margin40}
\end{align}
Then by \eqref{eq:q3-signed-fourier-integral} and \eqref{eq:q3-signed-margin40},
we know \eqref{eq:q3-direct-signed} holds.
\end{proof}

\begin{proof}[Proof of Theorem~\ref{thm:q3-main}]
Let $n\ge150$ be an integer, and let $p_n$ be the density of the sum
of $n$ independent uniform random variables on $[-1,1]$.
Fix $t>0$. By \eqref{eq:q3-aux40-values}, we have $w(0)>0$ and $w'(0)>0$.
Using \eqref{eq:q3-tilt-jetnumerator}, \eqref{eq:q3-aux80-score},
and Proposition~\ref{prop:q3-tilt-direct}, we obtain
\begin{align*}
 (\log w)'''(0)
 &=\frac{w'''(0)}{w(0)}
 +\frac{w'(0)}{w(0)}
 \left[-\frac{3w''(0)}{w(0)}+\frac{2w'(0)^2}{w(0)^2}\right]
 <\frac{w'''(0)+\frac72w'(0)}{w(0)}<0.
\end{align*}
By \eqref{eq:q3-tilt-log-third},
\[
 (\log p_n)'''(n\mathbb E[U_t])
 =\frac{(\log w)'''(0)}{[n\operatorname{Var}(U_t)]^{3/2}}<0.
\]
Since \eqref{eq:q3-tilted-mean-range} shows that $n\mathbb E[U_t]$
ranges over $(0,n)$ as $t$ ranges over $(0,\infty)$, for every $0<x<n$,
\[
 (\log p_n)'''(x)<0.
\]
For every $0<x<\sqrt n$, we have $0<\sqrt n\,x<n$.
By \eqref{eq:q3-density-rescaling},
\[
 (\log q_n)'''(x)
 =n^{\frac32}(\log p_n)'''(\sqrt n\,x)<0.
\]
Proposition~\ref{prop:cube-criterion} gives the entropy inequality
and its equality statement.
\end{proof}

\section{Applications to information theory}
\label{sec:information-theory}

Entropy bounds give estimates for the capacity of additive noise channels
and for the entropy power of sums. Madiman, Nayar, and Tkocz
\cite{MNT2021} and Melbourne, Nayar, and Roberto \cite[Section~2]{MNR2026}
used entropy lower bounds for log-concave variables to study these questions.
We use Theorem~\ref{thm:q3-main} to obtain capacity lower bounds that hold
for every choice of coefficients, compare the capacities in the high-power
limit, and determine the optimal reverse entropy power constant for a fixed
number of uniform summands.

Recall that, for a random variable $X$ with density $f$ and finite
differential entropy,
\begin{align*}
 h(X)&=-\mathbb E[\log f(X)]=-\int_{-\infty}^{\infty}f(x)\log f(x)\dd x.
\end{align*}
Its entropy power is
\begin{align*}
 N(X)&=\frac1{2\pi e}\exp[2h(X)].
\end{align*}
Thus $N(X)$ is the variance of a Gaussian variable with the same entropy as $X$.
Translation and scaling of the density give
$N(bX+c)=b^2N(X)$ for every real $b\ne0$ and $c$.
For independent random variables $X,Y$ with finite entropies and variances,
the entropy power inequality of Shannon and Stam \cite{Shannon1948,Stam1959} gives
\begin{align}\label{eq:info-epi}
 N(X+Y)&\ge N(X)+N(Y).
\end{align}

Throughout this section, let $n\ge150$ be an integer, and let
$U_1,\ldots,U_n$ be independent and uniform on $[-1,1]$.
For the equal-coefficient sum, whose density is $q_n$ from
Section~\ref{sec:cube-comparison}, let
\begin{align*}
 v_n&=N\left(\frac1{\sqrt n}\sum\limits_{i=1}^n U_i\right).
\end{align*}
For every $a\in\mathbb R^n$ with $\sum\limits_{i=1}^n a_i^2=1$,
the sum $\sum\limits_{i=1}^n a_iU_i$ has variance $1/3$.
Since $n^{-1/2}\sum\limits_{i=1}^n U_i$ is not Gaussian,
Theorem~\ref{thm:q3-main} and the Gaussian maximum entropy inequality
give
\begin{align}\label{eq:info-maximal-entropy-power}
 0<N\left(\sum\limits_{i=1}^n a_iU_i\right)&\le v_n<\frac13.
\end{align}

\subsection{Additive noise channels}

Let $Z$ have a bounded density $f_Z$ and finite variance, and let $P>0$.
An additive noise channel has input $X$ and output $X+Z$, where $X$
is independent of $Z$. We consider inputs with $\operatorname{Var}(X)\le P$.
Given $X=x$, the conditional output density is $y\mapsto f_Z(y-x)$.
Averaging over $X$ gives the marginal output density
\begin{align*}
 f_{X+Z}(y)&=\mathbb E[f_Z(y-X)]\le\|f_Z\|_\infty.
\end{align*}
Since $\operatorname{Var}(X+Z)\le P+\operatorname{Var}(Z)$, the Gaussian
maximum entropy inequality gives
\begin{align}\label{eq:info-output-entropy}
 -\log\|f_Z\|_\infty
 &\le h(X+Z)\le\frac12\log[2\pi e(P+\operatorname{Var}(Z))]<\infty.
\end{align}
The same argument with $X=0$ shows that $h(Z)$ is finite.

The mutual information is the expected log ratio of the conditional
output density to the marginal output density. Here it equals
\begin{align*}
 \mathbb E\!\left[\log\frac{f_Z(Z)}{f_{X+Z}(X+Z)}\right]
 &=h(X+Z)-h(Z).
\end{align*}
The capacity under the variance constraint $\operatorname{Var}(X)\le P$
is the supremum of this entropy difference:
\begin{align*}
 C_P(Z)&=\sup_{\substack{X\ \mathrm{independent\ of}\ Z\\
                         \operatorname{Var}(X)\le P}}
             [h(X+Z)-h(Z)].
\end{align*}
The supremum is over all input distributions, including discrete ones.
Translating the input leaves the entropy difference unchanged, so centering
$X$ shows that the variance constraint gives the same capacity as
$\mathbb E[X^2]\le P$. By \eqref{eq:info-output-entropy}, the capacity is finite.

We first recall Shannon's capacity bounds \cite[Theorem~18]{Shannon1948}
in the scalar formulation given in \cite[Section~2.1]{MNT2021}.
Let $G$ be centered Gaussian with variance $P$, independent of
$(Z,U_1,\ldots,U_n)$.
Since $N(G)=P$, \eqref{eq:info-epi} gives
\begin{align}\label{eq:info-capacity-lower}
 C_P(Z)&\ge h(G+Z)-h(Z)=\frac12\log\frac{N(G+Z)}{N(Z)}
 \ge\frac12\log\left(1+\frac P{N(Z)}\right).
\end{align}
By \eqref{eq:info-output-entropy} and the definition of $C_P(Z)$, we obtain
\begin{align}\label{eq:info-capacity-upper}
 C_P(Z)&\le\frac12\log[2\pi e(P+\operatorname{Var}(Z))]-h(Z)
 =\frac12\log\left(\frac{P+\operatorname{Var}(Z)}{N(Z)}\right).
\end{align}
Theorem~\ref{thm:q3-main} now gives a lower bound that holds for every
normalized coefficient vector. It also bounds the difference between the
capacity for equal coefficients and the infimum over all such vectors.

\begin{corollary}\label{cor:info-capacity}
Let $n\ge150$ be an integer. For every $P>0$ and every
$a\in\mathbb R^n$ with $\|a\|_2=1$,
\begin{align}\label{eq:info-uniform-capacity}
 C_P\left(\sum\limits_{i=1}^n a_iU_i\right)&\ge\frac12\log\left(1+\frac P{v_n}\right).
\end{align}
Moreover, for every $P>0$,
\begin{align}\label{eq:info-capacity-infimum}
 0&\le C_P\left(\frac1{\sqrt n}\sum\limits_{i=1}^n U_i\right)-\inf_{\|a\|_2=1}C_P\left(\sum\limits_{i=1}^n a_iU_i\right)
 \le\frac12\log\left(\frac{P+\frac13}{P+v_n}\right).
\end{align}
Uniformly over all such coefficient vectors $a$,
\begin{align}\label{eq:info-capacity-limit}
 \lim_{P\to\infty}\left[C_P\left(\sum\limits_{i=1}^n a_iU_i\right)-C_P\left(\frac1{\sqrt n}\sum\limits_{i=1}^n U_i\right)\right]&=h\left(\frac1{\sqrt n}\sum\limits_{i=1}^n U_i\right)-h\left(\sum\limits_{i=1}^n a_iU_i\right).
\end{align}
For each fixed coefficient vector $a$ with $\|a\|_2=1$ whose coordinates
are not all equal in magnitude, and for all sufficiently large $P$, we have
\begin{align*}
 C_P\left(\sum\limits_{i=1}^n a_iU_i\right)
 &>C_P\left(\frac1{\sqrt n}\sum\limits_{i=1}^n U_i\right).
\end{align*}
\end{corollary}

\begin{proof}
Fix $P>0$.
By \eqref{eq:info-capacity-lower} and \eqref{eq:info-maximal-entropy-power},
for every $a\in\mathbb R^n$ with $\|a\|_2=1$,
\begin{align}\label{eq:info-capacity-chain}
 \frac12\log\left(1+\frac P{v_n}\right)
 &\le h\left(G+\sum\limits_{i=1}^n a_iU_i\right)-h\left(\sum\limits_{i=1}^n a_iU_i\right)
 \le C_P\left(\sum\limits_{i=1}^n a_iU_i\right).
\end{align}
This proves \eqref{eq:info-uniform-capacity}.
The equal-coefficient vector is among those in the infimum in
\eqref{eq:info-capacity-infimum}, which proves nonnegativity in
\eqref{eq:info-capacity-infimum}.
Since the equal-coefficient noise has variance $1/3$ and entropy power
$v_n$, \eqref{eq:info-capacity-upper} and \eqref{eq:info-capacity-chain} give
\begin{align*}
 C_P\left(\frac1{\sqrt n}\sum\limits_{i=1}^n U_i\right)
 -\inf_{\|a\|_2=1}C_P\left(\sum\limits_{i=1}^n a_iU_i\right)
 &\le\frac12\log\frac{P+\frac13}{v_n}
 -\frac12\log\left(1+\frac P{v_n}\right)
 =\frac12\log\frac{P+\frac13}{P+v_n}.
\end{align*}
This proves \eqref{eq:info-capacity-infimum}.

Adding $h(Z)=\frac12\log[2\pi eN(Z)]$ to the bounds in
\eqref{eq:info-capacity-lower} and \eqref{eq:info-capacity-upper}, and using
$N(Z)>0$, we obtain
\begin{align*}
 \frac12\log(2\pi eP)
 &\le\frac12\log\!\bigl[2\pi e(P+N(Z))\bigr]
 \le C_P(Z)+h(Z)
 \le\frac12\log\!\bigl[2\pi e(P+\operatorname{Var}(Z))\bigr].
\end{align*}
The two noises in \eqref{eq:info-capacity-limit} both have variance $1/3$.
Applying the upper bound to one and the lower bound to the other, and
then interchanging them, gives
\begin{align}
 &\left|C_P\left(\sum\limits_{i=1}^n a_iU_i\right)
       +h\left(\sum\limits_{i=1}^n a_iU_i\right)
       -C_P\left(\frac1{\sqrt n}\sum\limits_{i=1}^n U_i\right)
       -h\left(\frac1{\sqrt n}\sum\limits_{i=1}^n U_i\right)\right|\notag\\
 \le{}&\frac12\log\left[2\pi e\left(P+\frac13\right)\right]
       -\frac12\log(2\pi eP)
 =\frac12\log\left(1+\frac1{3P}\right)\xrightarrow[P\to\infty]{}0.
 \label{eq:info-uniform-high-power}
\end{align}
Since this upper bound is independent of $a$, it proves the uniform
convergence in \eqref{eq:info-capacity-limit}.

Finally, fix a coefficient vector $a$ with $\|a\|_2=1$ whose coordinates
are not all equal in magnitude. By the equality statement in
Theorem~\ref{thm:q3-main},
\begin{align*}
 h\left(\frac1{\sqrt n}\sum\limits_{i=1}^n U_i\right)
 -h\left(\sum\limits_{i=1}^n a_iU_i\right)&>0.
\end{align*}
By \eqref{eq:info-uniform-high-power}, for all sufficiently large $P$ we have
\begin{align*}
 C_P\left(\sum\limits_{i=1}^n a_iU_i\right)
 -C_P\left(\frac1{\sqrt n}\sum\limits_{i=1}^n U_i\right)
 &\ge h\left(\frac1{\sqrt n}\sum\limits_{i=1}^n U_i\right)
       -h\left(\sum\limits_{i=1}^n a_iU_i\right)
       -\frac12\log\left(1+\frac1{3P}\right)>0.
\end{align*}
\end{proof}

\subsection{A sharp reverse entropy power inequality}

For independent uniform summands, Theorem~\ref{thm:q3-main} gives an
upper bound complementary to \eqref{eq:info-epi}.

\begin{corollary}\label{cor:info-reverse-epi}
Let $n\ge150$ be an integer, and let $X_1,\ldots,X_n$ be independent
random variables, each uniform on an interval of positive finite length. Then
\begin{align}\label{eq:info-reverse-epi}
 N\left(\sum\limits_{i=1}^n X_i\right)
 &\le\frac{1}{2}\pi e v_n\sum\limits_{i=1}^n N(X_i).
\end{align}
The constant $\pi e v_n/2$ is optimal. Equality holds if and only if
all $n$ intervals have the same length.
\end{corollary}

\begin{proof}
By translation invariance, we may take $X_i=b_iU_i$ with $b_i>0$.
Since $h(X_i)=\log(2b_i)$, we have $N(X_i)=2b_i^2/(\pi e)$.
By scaling and \eqref{eq:info-maximal-entropy-power},
\begin{align*}
 \frac{N\left(\sum\limits_{i=1}^n X_i\right)}{\sum\limits_{i=1}^n N(X_i)}
 &=\frac{1}{2}\pi e\,
   \frac{N\left(\sum\limits_{i=1}^n b_iU_i\right)}{\sum\limits_{i=1}^n b_i^2}
 =\frac{1}{2}\pi e\,
   N\left(\frac{\sum\limits_{i=1}^n b_iU_i}{\left(\sum\limits_{i=1}^n b_i^2\right)^{1/2}}\right)
 \le\frac{1}{2}\pi e v_n.
\end{align*}
By Theorem~\ref{thm:q3-main}, equality holds exactly when all $b_i$ are equal,
or equivalently, when all intervals have the same length.
These choices also show that the constant is optimal.
\end{proof}

The Gaussian maximum entropy inequality gives
\begin{align*}
 N\left(\sum\limits_{i=1}^n X_i\right)
 &\le\operatorname{Var}\left(\sum\limits_{i=1}^n X_i\right)
 =\sum\limits_{i=1}^n\operatorname{Var}(X_i)
 =\frac{\pi e}{6}\sum\limits_{i=1}^n N(X_i).
\end{align*}
The same constant $\pi e/6$ appears in the reverse entropy power
inequality for uncorrelated symmetric log-concave variables
\cite[Corollary~18]{MNT2021}.
Since $v_n<1/3$, the optimal constant in \eqref{eq:info-reverse-epi}
is strictly smaller for every fixed $n\ge150$.

\appendix

\section{Moment inequalities}
\label{sec:q3-even-moments80}

For the tilted variable $Y$ of Section~\ref{sec:q3-tilt-reduction}, we prove
all three bounds in Lemma~\ref{lem:q3-moments}.

\begin{proof}[Proof of Lemma~\ref{lem:q3-moments}]
Fix $t>0$. The variable $Y$ is bounded, and its density satisfies
\eqref{eq:q3-tilted-density-sign}, with a strictly positive difference on an
interval. Since $y\mapsto y^3$ and $y\mapsto y^5$ are odd and strictly convex
on the positive half-line, \eqref{eq:q3-odd-convex-comparison} gives
$\mathbb E[Y^3]<0$ and $\mathbb E[Y^5]<0$, and thus
$\mathbb E[Y^5]/\mathbb E[Y^3]>0$.
Since $K'''(t)=[K''(t)]^{3/2}\mathbb E[Y^3]<0$, the variance $K''(t)$
decreases with $t$. The variable $U_t$ is not constant, so
$0<\operatorname{Var}(U_t)=K''(t)\le K''(0)=1/3$.
The variable $-Y$ has a log-concave density, mean zero, and unit variance.
By \cite[Lemma~2]{BubeckEldan2015},
$-\mathbb E[Y^3]=\mathbb E[(-Y)^3]\le2$.

For the fifth-moment ratio, \eqref{eq:q3-tilted-mgf} and
$\mathbb E[U_t]=K'(t)$ give, as $s\to0$,
\begin{align*}
 &\mathbb E\left[\exp\left[s(U_t-\mathbb E[U_t])\right]\right]\\
 ={}&\exp[K(t+s)-K(t)-sK'(t)]\\
 ={}&\exp\left[\frac{K''(t)}{2!}s^2+\frac{K'''(t)}{3!}s^3
 +\frac{K^{(4)}(t)}{4!}s^4+\frac{K^{(5)}(t)}{5!}s^5+O(s^6)\right].
\end{align*}
Comparing the coefficients of $s^5$ gives
\[
 \frac{\mathbb E[(U_t-\mathbb E[U_t])^5]}{5!}
 =\frac{K^{(5)}(t)}{5!}+\frac{K''(t)K'''(t)}{2!\,3!}.
\]
Dividing the fifth centered moment by $[K''(t)]^{5/2}$ and using
$\mathbb E[Y^3]=K'''(t)/[K''(t)]^{3/2}$, we obtain
\begin{align*}
 \frac{\mathbb E[Y^5]}{\mathbb E[Y^3]}
 &=\frac{K^{(5)}(t)+10K''(t)K'''(t)}{K''(t)K'''(t)}
 =10+\frac{K^{(5)}(t)}{K''(t)K'''(t)}.
\end{align*}
Since $K''(t)>0$ and $K'''(t)<0$, the denominator is negative, so
\begin{equation}\label{eq:q3-fifth-moment-criterion}
 \frac{\mathbb E[Y^5]}{\mathbb E[Y^3]}<22
 \quad\Longleftrightarrow\quad
 K^{(5)}(t)-12K''(t)K'''(t)>0.
\end{equation}
For $t>0$, differentiation of \eqref{eq:q3-log-mgf} gives
\begin{align}
 K''(t)&=\frac1{t^2}-\frac1{\sinh^2t},\quad
 K'''(t)=-\frac2{t^3}+\frac{2\cosh t}{\sinh^3t},\quad
 K^{(5)}(t)=-\frac{24}{t^5}
 +\frac{8\cosh t}{\sinh^3t}+\frac{24\cosh t}{\sinh^5t}.
 \label{eq:q3-moment-derivatives}
\end{align}
By \eqref{eq:q3-moment-derivatives}, we obtain
\begin{align}
 &t^3\sinh^5t[K^{(5)}(t)-12K''(t)K'''(t)]\notag\\
 ={}&t^3\sinh^5t
 \biggl[8\left(1-\frac3{t^2}\right)\frac{\cosh t}{\sinh^3t}
 +\frac{48\cosh t}{\sinh^5t}-\frac{24}{t^3\sinh^2t}\biggr]\notag\\
 ={}&8t(t^2-3)\cosh t\,\sinh^2t+48t^3\cosh t-24\sinh^3t\notag\\
 ={}&2t(t^2-3)\cosh(3t)+2t(23t^2+3)\cosh t-6\sinh(3t)+18\sinh t\notag\\
 ={}&8\sum\limits_{k=4}^\infty
 \left[9^{k-1}(2k^3-14k-27)+46k^3-10k+3\right]
 \frac{t^{2k+1}}{(2k+1)!}\notag\\
 ={}&8\sum\limits_{k=4}^\infty
 \Bigl(9^{k-1}\bigl[45+2(k-4)(k^2+4k+9)\bigr]+k(46k^2-10)+3\Bigr)
 \frac{t^{2k+1}}{(2k+1)!}>0,
 \label{eq:q3odd:fifth-hyperbolic}
\end{align}
where we used the identities $4\cosh t\,\sinh^2t=\cosh(3t)-\cosh t$ and
$4\sinh^3t=\sinh(3t)-3\sinh t$, together with the Taylor series for
$\sinh$ and $\cosh$.
Since $t^3\sinh^5t>0$, \eqref{eq:q3odd:fifth-hyperbolic} gives
$K^{(5)}(t)-12K''(t)K'''(t)>0$, and
\eqref{eq:q3-fifth-moment-criterion} proves the required ratio bound.
\end{proof}

\section{Moment integrals}
\label{sec:q3-moment-tails80}

We prove the moment and integral bounds used in
Sections~\ref{sec:q3-tilt-even} and~\ref{sec:q3-tilt-direct}.
The moment integrals $M_{2k}(n)$ are defined in
\eqref{eq:q3-comparison-moments}, with
$r_n(u)=(1+u^2/n)^{-(n-1)/2}$.

\begin{lemma}\label{lem:q3-integrals80}
For every $n\ge150$,
\begin{equation}\label{eq:q3-aux80-moment-bounds}
 M_2(n)<\frac{21}{20},\qquad M_4(n)<5,\qquad
 M_6(n)<24,\qquad M_{10}(n)<1500,\qquad M_{12}(n)<16000.
\end{equation}
\end{lemma}

\begin{proof}
By \eqref{eq:q3-comparison-moments}, for $n>3$ we have
\begin{align}
 M_0(n)&=\left(\frac1{2\pi}\int_{-\infty}^{\infty}\int_{-\infty}^{\infty}
 r_n(x)r_n(y)\,\mathrm dx\,\mathrm dy\right)^{1/2}\notag\\
 &\le\left(\frac1{2\pi}\int_{-\infty}^{\infty}\int_{-\infty}^{\infty}
 \left(1+\frac{x^2+y^2}{n}\right)^{-\frac{n-1}{2}}
 \,\mathrm dx\,\mathrm dy\right)^{1/2}\notag\\
 &=\left(\int_0^\infty s\left(1+\frac{s^2}{n}\right)^{-\frac{n-1}{2}}
 \,\mathrm ds\right)^{1/2}\notag\\
 &=\left(\left[-\frac n{n-3}
 \left(1+\frac{s^2}{n}\right)^{-\frac{n-3}{2}}\right]_0^\infty\right)^{1/2}\notag\\
 &=\sqrt{\frac n{n-3}}=\sqrt{1+\frac3{n-3}}
 \le1+\frac3{2(n-3)}=\frac{2n-3}{2(n-3)},\label{eq:q3-aux80-zeroth-moments}
\end{align}
where the first inequality follows from
$(1+x^2/n)(1+y^2/n)\ge1+(x^2+y^2)/n$, the second equality uses a change
to polar coordinates, and the last inequality follows from concavity
of the square root.

For every integer $k\ge1$ and $n>2k+2$, we have
\begin{align*}
 M_{2k}(n)&=\frac1{\sqrt{2\pi}}\int_{-\infty}^{\infty}
 u^{2k}\left(1+\frac{u^2}{n}\right)^{-\frac{n-1}{2}}\,\mathrm du\\
 &=\frac{(2k-1)n}{(n-3)\sqrt{2\pi}}\int_{-\infty}^{\infty}
 u^{2k-2}\left(1+\frac{u^2}{n}\right)^{-\frac{n-3}{2}}\,\mathrm du\\
 &=\frac{(2k-1)n}{(n-3)\sqrt{2\pi}}\int_{-\infty}^{\infty}
 u^{2k-2}\left(1+\frac{u^2}{n}\right)r_n(u)\,\mathrm du\\
 &=\frac{(2k-1)n}{n-3}\left(M_{2k-2}(n)+\frac1nM_{2k}(n)\right),
\end{align*}
where the second equality follows from integration by parts, with
$u^{2k-1}(1+u^2/n)^{-(n-3)/2}\to0$ as $|u|\to\infty$ because $n>2k+2$.
Rearranging gives
\begin{align}\label{eq:q3-aux80-moment-recurrences}
 M_{2k}(n)&=\frac{(2k-1)n}{n-2k-2}M_{2k-2}(n).
\end{align}

For $n\ge150$, \eqref{eq:q3-aux80-zeroth-moments} and
\eqref{eq:q3-aux80-moment-recurrences} give
\[
 M_2(n)=\frac n{n-4}M_0(n)
 \le\frac n{n-4}\frac{2n-3}{2(n-3)}
 \le\frac{150}{146}\cdot\frac{99}{98}<\frac{21}{20}.
\]
Repeated use of \eqref{eq:q3-aux80-moment-recurrences}, followed by
\eqref{eq:q3-aux80-zeroth-moments}, gives
\begin{align*}
 M_{12}(n)&=M_0(n)\prod\limits_{j=1}^{6}\frac{(2j-1)n}{n-2j-2}
 \le10395\cdot\frac{2n-3}{2(n-3)}
 \prod\limits_{j=1}^{6}\frac n{n-2j-2}<16000.
\end{align*}
Since $(2k-1)n/(n-2k-2)>2k-1$, the same recurrence gives
\begin{align*}
 M_{10}(n)&<\frac{M_{12}(n)}{11}<1500,\qquad
 M_6(n)<\frac{M_{12}(n)}{11\cdot9\cdot7}<24,\qquad
 M_4(n)<\frac{M_{12}(n)}{11\cdot9\cdot7\cdot5}<5.
\end{align*}
This proves \eqref{eq:q3-aux80-moment-bounds}.
\end{proof}

The next lemma bounds the tail integrals, the contribution from
$[0,\sqrt{14}/2]$, and the error term in \eqref{eq:q3-central-loss200}.

\begin{lemma}\label{lem:q3-integration-errors}
Let $n\ge150$.

\noindent\textup{(i)} The tail integrals satisfy
\begin{equation}\label{eq:q3-tail-zero}
 \frac1{\sqrt{2\pi}}\int_{|u|>\frac{\sqrt n}{2}}
 \left[\exp\left(-\frac12u^2\right)+r_n(u)\right]\,\mathrm du<\frac1n,
\end{equation}
\begin{equation}\label{eq:q3-tail-four}
 \frac1{\sqrt{2\pi}}\int_{|u|>\frac{\sqrt n}{2}}
 u^4\left[\exp\left(-\frac12u^2\right)+r_n(u)\right]\,\mathrm du<\frac3n,
\end{equation}
\begin{equation}\label{eq:q3-tail-six}
 \frac1{\sqrt{2\pi}}\int_{|u|>\frac{\sqrt n}{2}}
 u^6\left[\exp\left(-\frac12u^2\right)+r_n(u)\right]\,\mathrm du<\frac35.
\end{equation}
\noindent\textup{(ii)} On $0\le u\le\sqrt{14}/2$, we have
\begin{equation}\label{eq:q3-envelope-near-zero}
 r_n(u)\le\frac{21}{20}\exp\left(-\frac12u^2\right),
\end{equation}
and
\begin{equation}\label{eq:q3-negative-weight-loss200}
 \sqrt{\frac2\pi}\int_0^{\frac{\sqrt{14}}{2}}u\left(\frac72-u^2\right)
 \frac{u^3}{6}\exp\left(-\frac12u^2\right)\,\mathrm du<1.
\end{equation}

\noindent\textup{(iii)} For the error term in \eqref{eq:q3-central-loss200},
we have
\begin{equation}\label{eq:q3-integrated-central-error}
 \sqrt{\frac2\pi}\int_{\frac{\sqrt{14}}{2}}^\infty u\left(u^2-\frac72\right)
 \left[\left(\frac{11u^5}{60}+\frac{u^7}{42}\right)
 \exp\left(-\frac12u^2\right)+\frac{u^9}{216}r_n(u)\right]\,\mathrm du<80.
\end{equation}
\end{lemma}

\begin{proof}
\textup{(i)} For every integer $k\ge0$, integration by parts gives
\begin{equation}\label{eq:q3-gaussian-moments}
 \frac1{\sqrt{2\pi}}\int_{-\infty}^{\infty}
 u^{2k}\exp\left(-\frac12u^2\right)\,\mathrm du
 =\frac{(2k)!}{2^k k!}.
\end{equation}
For $|u|>\sqrt n/2$, we have $u^2>n/4$.
Using \eqref{eq:q3-comparison-moments},
\eqref{eq:q3-aux80-moment-bounds}, and
\eqref{eq:q3-gaussian-moments}, we obtain
\begin{align*}
 &\frac1{\sqrt{2\pi}}\int_{|u|>\frac{\sqrt n}{2}}
 \left[\exp\left(-\frac12u^2\right)+r_n(u)\right]\,\mathrm du\\
 \le{}&\frac{64}{n^3\sqrt{2\pi}}\int_{-\infty}^{\infty}
 u^6\left[\exp\left(-\frac12u^2\right)+r_n(u)\right]\,\mathrm du
 =\frac{64}{n^3}\bigl(15+M_6(n)\bigr)<\frac1n,\\[6pt]
 &\frac1{\sqrt{2\pi}}\int_{|u|>\frac{\sqrt n}{2}}
 u^4\left[\exp\left(-\frac12u^2\right)+r_n(u)\right]\,\mathrm du\\
 \le{}&\frac{256}{n^4\sqrt{2\pi}}\int_{-\infty}^{\infty}
 u^{12}\left[\exp\left(-\frac12u^2\right)+r_n(u)\right]\,\mathrm du
 =\frac{256}{n^4}\bigl(10395+M_{12}(n)\bigr)<\frac3n,\\[6pt]
 &\frac1{\sqrt{2\pi}}\int_{|u|>\frac{\sqrt n}{2}}
 u^6\left[\exp\left(-\frac12u^2\right)+r_n(u)\right]\,\mathrm du\\
 \le{}&\frac{64}{n^3\sqrt{2\pi}}\int_{-\infty}^{\infty}
 u^{12}\left[\exp\left(-\frac12u^2\right)+r_n(u)\right]\,\mathrm du
 =\frac{64}{n^3}\bigl(10395+M_{12}(n)\bigr)<\frac35.
\end{align*}

\textup{(ii)} For $0\le u\le\sqrt{14}/2$, the inequality
$\log(1+x)\ge x-x^2/2$ for $x\ge0$ gives
\begin{align*}
 r_n(u)\exp\left(\frac12u^2\right)
 &=\exp\left[\frac{u^2}{2}-\frac{n-1}{2}\log\left(1+\frac{u^2}{n}\right)\right]
 \le\exp\left(\frac{u^2}{2n}+\frac{u^4}{4n}\right)<\frac{21}{20},
\end{align*}
which proves \eqref{eq:q3-envelope-near-zero}.
\begin{align*}
 \sqrt{\frac2\pi}\int_0^{\frac{\sqrt{14}}{2}}u\left(\frac72-u^2\right)
 \frac{u^3}{6}\exp\left(-\frac12u^2\right)\,\mathrm du
 &\le\frac16\sqrt{\frac2\pi}\int_0^{\frac{\sqrt{14}}{2}}
 \left(\frac72u^4-u^6\right)\,\mathrm du
 =\frac{49}{120}\sqrt{\frac7\pi}<1.
\end{align*}
This proves \eqref{eq:q3-negative-weight-loss200}.

\textup{(iii)} By evenness and \eqref{eq:q3-gaussian-moments},
\begin{align*}
 \sqrt{\frac2\pi}\int_0^\infty u\left(u^2-\frac72\right)
 \left(\frac{11u^5}{60}+\frac{u^7}{42}\right)
 \exp\left(-\frac12u^2\right)\,\mathrm du&=\frac{187}{8}.
\end{align*}
By \eqref{eq:q3-comparison-moments}, \eqref{eq:q3-aux80-moment-bounds},
and \eqref{eq:q3-aux80-moment-recurrences}, we obtain
{
\medmuskip=3mu\relax
\thickmuskip=4mu\relax
\begin{align*}
 \sqrt{\frac2\pi}\int_0^\infty
 u\left(u^2-\frac72\right)\frac{u^9}{216}r_n(u)\,\mathrm du
 &=\frac1{216}\left(M_{12}(n)-\frac72M_{10}(n)\right)
 =\frac{M_{12}(n)}{216}\left[1-\frac{7(n-14)}{22n}\right]<53.
\end{align*}
}
Since both $r_n(u)$ and $\exp(-u^2/2)$ are at most one, direct integration gives
\begin{align*}
 &-\sqrt{\frac2\pi}\int_0^{\frac{\sqrt{14}}{2}}u\left(u^2-\frac72\right)
 \left[\left(\frac{11u^5}{60}+\frac{u^7}{42}\right)
 \exp\left(-\frac12u^2\right)+\frac{u^9}{216}r_n(u)\right]\,\mathrm du\\
 \le{}&-\sqrt{\frac2\pi}\int_0^{\frac{\sqrt{14}}{2}}u\left(u^2-\frac72\right)
 \left(\frac{11u^5}{60}+\frac{u^7}{42}+\frac{u^9}{216}\right)\,\mathrm du<3.
\end{align*}
Summing these three bounds proves \eqref{eq:q3-integrated-central-error}.
\end{proof}

\end{document}